\documentclass[12pt,draftcls,onecolumn]{IEEEtran}

\usepackage{amssymb,latexsym,amsfonts,amsmath}
\usepackage{amsthm}
\usepackage{multirow} 
\usepackage{graphicx,diagrams}
\usepackage{comment}
\usepackage{epsfig}
\usepackage{graphicx}
\usepackage{color}
\graphicspath{figures}

\newtheorem{theorem}{Theorem}[section]
\newtheorem{lemma}[theorem]{Lemma}

\newtheorem{proposition}[theorem]{Proposition}
\newtheorem{corollary}[theorem]{Corollary}

\newtheorem{definition}[theorem]{Definition}
\newtheorem{example}[theorem]{Example}

\newtheorem{remark}[theorem]{Remark}
\numberwithin{equation}{section}

\renewcommand{\Re}{{\mathbb{R}}}
\newcommand{\R}{{\mathbb{R}}}

\newcommand{\N}{{\mathbb{N}}}

\newcommand{\ie}{{\it i.e.}}
\newcommand{\eg}{{\it e.g. }}

\newcommand{\s}{\mathrm{span}}

\newcommand{\eq}[1]{\begin{equation}#1\end{equation}}
\newcommand{\eqs}[1]{\begin{eqnarray}#1\end{eqnarray}}
\renewcommand{\exp}{\mathbf{e}}
\newcommand{\e}{\mathrm{e}}
\newcommand{\y}{\mathrm{y}}
\renewcommand{\u}{\mathrm{u}}
\newcommand{\x}{\mathrm{x}}
\newcommand{\w}{\mathrm{w}}

\newcommand{\z}{\mathrm{z}}
\renewcommand{\s}{\mathrm{s}}
\newcommand{\diffV}{\frac{\partial V}{\partial x}}

\title{Exploiting isochrony in self-triggered control}

\author{\parbox{4 in}{\centering Adolfo Anta \qquad Paulo Tabuada\thanks{This research was partially supported by the National Science Foundation award 0834771 and the Alexander von Humboldt Foundation.}
\thanks{P. Tabuada is with the Department of Electrical Engineering, University of California at Los Angeles, Los Angeles, CA 90095-1594, USA: {\tt\small tabuada@ee.ucla.edu}}%
\thanks{A. Anta is with the Technische Universit\"at Berlin \& Max Planck Institute f\"ur Dynamik komplexer technischer Systeme, Germany: {\tt\small anta@control.tu-berlin.de}}}}
\begin{document}
\maketitle

\begin{abstract}
Event-triggered control and self-triggered control have been recently proposed as new implementation paradigms 
that reduce resource usage for control systems. In self-triggered control, the controller is augmented with the computation of the next time instant at which the feedback control law is to be recomputed. Since these execution instants are obtained as a function of the plant state, we effectively close the loop only when it is required to maintain the desired performance, thereby greatly reducing the resources required for control. In this paper we present a new technique for the computation of the execution instants by exploiting the concept of isochronous manifolds, also introduced in this paper. While our previous results showed how homogeneity can be used to compute the execution instants along \emph{some} directions in the state space, the concept of isochrony allows us to compute the executions instants along \emph{every} direction in the state space. Moreover, we also show in this paper how to homogenize smooth control systems thus making our results applicable to any smooth control system. The benefits of the proposed approach with respect to existing techniques are analyzed in two examples.
\end{abstract}

\section{Introduction}
Recent technological and economical trends require control systems not to be implemented on dedicated platforms, but rather on shared devices where different subsystems compete for the available resources (\eg computation time, communication bandwidth, shared actuators). In this shared setup, traditional implementation paradigms are no longer adequate since they require the same amount of resources independently of the workload, performance, or the effect of exogenous disturbances. For instance, in a digital implementation, control laws are executed at the same rate regardless of the current state of the system. While this approach simplifies the analysis and implementation of the control system, the choice of the period for the control loop is independent of the current state and performance of the system. To overcome the shortcomings of periodic implementations, several researchers advocated the use of event-triggered control~\cite{arzen99,astrom2002cra,tabuada07,heemels08,lunze}. Under this paradigm, control tasks are executed whenever a triggering condition is satisfied (\eg whenever the state deviates from a nominal trajectory by a certain amount). The use of event-triggered control, however, entails one handicap: dedicated hardware is required to continuously check the triggering condition. Even if such hardware can be deployed, it implies a significant cost. To surmount this drawback, Velasco et al~\cite{velasco2003stt} suggested the self-triggered control concept, by which the next execution time is given by a function of 	the last measurement of the state. Self-triggered control can be regarded as a software-based emulation of event-triggered control. Self-triggered control for linear systems was studied in~\cite{wang2008_tac,mazo2010aut}. It has also been recently applied in the context of sensor/actuator networks~\cite{manolo,tiberi-adaptive,camacho-self}. To the best of our knowledge, self-triggered formulas for nonlinear systems were first developed in~\cite{tabuada10_tac}.


In this paper we extend our previous work on self-triggered control in several directions. In~\cite{tabuada10_tac}, the geometric properties of homogeneous vector fields were exploited to derive 
scaling laws describing the inter-execution times, that is, the times elapsed between consecutive executions of the control law. These results explained how such times evolve \textit{along} homogeneous rays. In this paper we develop formulas that explain how the inter-execution times evolve \textit{across} homogeneous rays. The contributions of this paper are as follows:

\begin{itemize}

\item We introduce the concept of isochronous manifolds, show their existence, and how it can be used to describe how the inter-execution times change across homogeneous rays.


\item We combine our previous results in~\cite{tabuada10_tac} with the concept of isochronous manifolds to obtain new self-triggering formulas that explain how the inter-execution times evolve along rays and across rays.

\item We generalize the results in~\cite{tabuada10_tac} to non-homogeneous control systems and non-homogeneous triggering conditions.

\end{itemize}
We also analyze the accuracy of the proposed self-triggering formulas by:

\begin{itemize}
\item Providing a bound for the mismatch between the inter-execution times generated by the event-triggering technique and by the self-triggering technique. 
\item Analyzing the inherent trade-off between complexity and accuracy of the self-triggering formula, and deriving an iterative algorithm that trades complexity for precision.
\item Discussing the role played by the different design parameters in the accuracy of the self-triggering formulas.
\end{itemize}

This paper does not focus on the self-triggered emulation of a particular event-triggered implementation, but rather develops a self-triggering technique that can be used with any \textit{smooth event-triggering condition} defined for any closed-loop system described by a \textit{smooth vector field}. We illustrate the results on two examples, showing that the proposed technique outperforms previous solutions and produces self-triggered times that hardly differ from the event-triggered times. A preliminary version of these results appeared in~\cite{anta_isochronous}.

\section{Notation}
We use the notation $\vert x\vert$ to denote the Euclidean norm of an element $x\in\Re^n$. 
A function is said to be of class $C^{\infty}$ or smooth if it can be differentiated infinitely many times. All the objects in this paper are considered to be smooth unless otherwise stated.

Given vector fields $X$ and $Y$ in an n-dimensional manifold $M$, we let $[X,Y]$ denote their Lie product which, in local coordinates $x = (x_1,x_2,...,x_n)$, we take as $\frac{\partial{Y}}{\partial{x}}(x)X(x) - \frac{\partial{X}}{\partial{x}}(x)Y(x)$. We use ${\cal L}_X h(x)$ to denote the Lie derivative of a map $h:M\rightarrow \R$ at a point $x$ along the flow of the vector field $X:M \rightarrow TM$ which, in local coordinates, we take as $\frac{\partial{h}}{\partial{x}}X(x)$. Likewise, ${\cal L}_X^k h$ represents the $k$th Lie derivative, defined by ${\cal L}_X^0 h = h$ and ${\cal L}_X^k h = {\cal L}_X ({\cal L}_X^{k-1} h)$. 
Let $\phi:M \rightarrow N$ be a map, and let $X$ and $Y$ be vector fields in manifolds $M$ and $N$, respectively. We denote the differential of a map $\phi$ by $T\phi$, which, in local coordinates, we take as $\frac{\partial{\phi}}{\partial{x}}$. We call $X$ and $Y$ $\phi$-related if the following holds:
\eq{ 
\label{phi_rel_condition}
T\phi \cdot X = Y \circ \phi.}
In this paper we consider control systems of the form:
\begin{equation}
\label{CS}
\dot{\mathrm{x}}=f(\mathrm{x},\mathrm{u}),\qquad \mathrm{x}(t)\in \Re^n,\,\, \mathrm{u}(t)\in \R^m, t\in\R_0^+.
\end{equation}
We denote by $x \in \R^n$ the state of the control system, by $\mathrm{x}$ a solution of~(\ref{CS}), and by $\mathrm{u}$ the input trajectory. 
Finally, we use $\exp^A$ to represent the exponential of a matrix $A\in\R^{n\times n}$, $A_{ij}$ to denote the element $ij$ of the matrix $A$ and $I$ to represent the identity matrix.\\

\section{Self-triggering stabilization of nonlinear systems}
We start by analyzing the behavior of a control system under a sample-and-hold implementation. Consider the system~(\ref{CS}), for which a feedback control law of the form $u=k(x)$ has been designed. The implementation of such control law on an embedded processor is typically done by sampling the current state of the system $\x(t)$ at time instants $t=t_1, t_2, \ldots$,  computing the new value for the input $\u(t_i) = k(\x(t_i))$, and updating the actuator with $\u(t_i)$ at time instants $t_i$. We shall use the term \textit{execution} to denote this process of sampling the state, computing the control law and updating the input. Between actuator updates the input $u$ is held constant, according to:
\begin{equation}
t\in[t_i, t_{i+1}[ \Rightarrow \mathrm{u}(t)=\mathrm{u}(t_i).
\end{equation}
Traditionally, the controller is executed periodically, that is $t_{i+1}-t_i = T$, for any $i\in \N$, where $T>0$ is the period. Instead, we would like to identify a sequence of execution time instants that guarantees stability and desired performance while reducing the number of executions. 
To proceed, we define the measurement error $\e$ as the difference between the last measurement of the state and the current value of the state:
\begin{equation}
\label{Error}
\mathrm{e}(t)=\mathrm{x}(t_i)-\mathrm{x}(t) \qquad \text{for } t\in [t_i, t_{i+1}[.
\end{equation}
With this definition, the closed loop $\dot{\mathrm{x}} = f(\mathrm{x},k(\mathrm{x}(t_i)))$ becomes:
\begin{equation}
\label{closed_loop_dig}
\dot{\mathrm{x}}=f(\mathrm{x},k(\mathrm{x}+\mathrm{e})).
\end{equation}
Whenever a control law $u=k(x)$ is designed to render the system~(\ref{CS}) asymptotically stable, there exist a Lyapunov function $V$ satisfying:
\eq{
\label{ineq_V}
\dot V\vert_{t=0} = \frac{\partial V}{\partial x}f(x,k(x))<0 \quad \forall x\in\R^n \backslash \{0\}.
}
Moreover, control laws are usually designed not only to guarantee stability but also to achieve a certain performance. In this paper the performance requirement will be given by a desired rate of decay \mbox{$h:\R^n\rightarrow\R$} of the Lyapunov function. Hence $V$ satisfies:
\begin{equation}
\label{spec}
\diffV f(x,k(x)) \leq h(x)<0, \quad \forall x\in\R^n \backslash \{0\}.
\end{equation}
However, under a digital implementation the dynamics of the system are given by~(\ref{closed_loop_dig}). Let the desired rate of decay for the implemented closed loop be $\tilde{h}(x)>h(x)$. If the error $e$ satisfies:
\eq{
\label{ineq_V_e}
\frac{\partial V}{\partial x}f(x,k(x+e))\leq \tilde{h}(x)<0
}
the rate of decay $\tilde h$ can be guaranteed for the sample-and-hold implementation. Therefore, the sequence of time instants $t_i$ has to be chosen so that inequality~(\ref{ineq_V_e}) is satisfied at any time. This can be achieved by closing the loop whenever:
\eq{
\label{trig_condition}
\Gamma(x,e) := \frac{\partial V}{\partial x}f(x,k(x+e))-\tilde{h}(x) =0.
}
Equation~(\ref{trig_condition}) implicitly defines the time instants $\tau$ at which the state $x$ needs to be measured, the control law $u=k(x)$ computed, and the actuator updated with the control input $u$:
$$\tau(\x(t_i),\e(t_i)) := \min\{t>t_{i}:\Gamma(\x(t,\x(t_i),\e(t_i)),\e(t,\x(t_i),\e(t_i))) =0\}.$$
Whenever there does not exist $t$ such that $\Gamma(\x(t,\x(t_i),\e(t_i)),\e(t,\x(t_i),\e(t_i))) =0$, we will say that $\tau(\x(t_i),\e(t_i)) = \infty$. Indeed, if (\ref{trig_condition}) is never satisfied the control law does not need to be updated anymore. Notice that our analysis is not constrained by the particular structure of~\eqref{trig_condition} but extends to other typical (smooth) triggering conditions available in the literature, as it will be shown in the remainder of the paper.

Equation~(\ref{trig_condition}) represents an \textit{event-triggering condition} that renders the system asymptotically stable and guarantees the desired rate of decay $\tilde h$.
Upon the execution of the control law at $t=t_i$, the state is measured and the error becomes 0, since \mbox{$\mathrm{x}(t)=\mathrm{x}(t_i)$} implies \mbox{$\mathrm{e}(t_i)=\mathrm{x}(t_i)-\mathrm{x}(t_i)=0$}. 
An event-triggered implementation based on this equality would require testing~(\ref{ineq_V_e}) frequently. Unless this testing process is implemented in hardware, one might run the risk of consuming the processor time testing~(\ref{ineq_V_e}). To overcome this drawback, we consider self-triggering strategies, where the current state measurement $\mathrm{x}(t_i)$ is used to determine the next execution time for the control law. Self-triggered implementations represent a software-based emulation of event-triggered implementations.

For clarity of exposition, at this point we define a new state $z\in\R^{2n}$ that includes the state of the dynamical system $x$ and the measurement error $e$:
$$z := \left[\begin{array} {cc} x \\ e\end{array}\right].$$ 
Likewise, we use $\z$ to denote the trajectories of the extended state $z$.

The self-triggering strategy decides the next inter-execution time according to a function \mbox{$\tau^{\downarrow}:\R^{2n}\rightarrow \R^+$} that lower bounds the event-triggered times, that is, the inter-execution times of the event-triggering policy:
\eq{
\tau^{\downarrow}(\z(t_i)) \leq \tau(\z(t_i))
}
where $\z(t_i)=(\x(t_i),\e(t_i))$ corresponds to the last measurement of the current state $\x(t_i)$ and to the error at the sampling time $\e(t_i)$. Since
\eq{
\Gamma(\z(t,\z(t_i))) \leq 0 \qquad \forall t\in [t_i, \tau[
}
closing the loop at $\tau$, or before, enforces~\eqref{ineq_V_e} for all $t\geq 0$.

Therefore, the self-triggered times should be as large as possible (in order to reduce the number of executions of the control law) but no larger than the event-triggered times (in order to guarantee the desired performance). Since the state trajectories are in general not known for nonlinear systems, it is not possible to compute the inter-execution times $\tau(\z(t_i))$ in closed form. In the next section we develop scaling laws for the inter-execution times of smooth nonlinear systems implicitly defined by the triggering condition~\eqref{trig_condition}. 

\section{Scaling laws for the inter-execution times}
\label{sec_scal_laws}
We start by reviewing the scaling laws previously derived in~\cite{tabuada10_tac} for the special class of homogeneous systems.
\subsection{Homogeneous systems}
\label{subsec_scal_laws}
We first define homogeneity for functions.
\begin{definition}
\label{def_hom_fcn}
A function $f: \R^n \rightarrow \R^n$ is called homogeneous of degree $\xi\in\R$ if there exist $r_i > 0,\: i = 1\ldots n$ such that:
\begin{eqnarray}
\label{Dil_out}
\forall\lambda > 0 \qquad f_i (\lambda^{r_1} x_1,\ldots,\lambda^{r_n} x_n) = \lambda^{\xi} \lambda^{r_i} f_i(x_1,\ldots,x_n)
\end{eqnarray}
where $\xi > -\min_i{r_i}$. 
\end{definition}
In a geometric setting, a function $f$ is said to be homogeneous of degree $\xi$ if the following holds:
\eq{
\label{def_hom_geom}
{\cal L}_D f = (\xi+1)f
}
where $D$ is known as a dilation vector field, \ie, a vector field such that $\dot{\x} = - D(\x)$ is globally asymptotically stable. When the choice
\begin{equation}
\label{gen_dilation_vf}
D = \sum_{i=1}^n x_i \frac{\partial}{\partial x_i}
\end{equation}
is made, the geometric definition in~\eqref{def_hom_geom} results in Definition~\ref{def_hom_fcn} (see~\cite{kawski95,arnol'd_ode} for details). Likewise, we define the notion of homogeneity for vector fields.
\begin{definition}~\cite{kawski95}
Let $D:M \rightarrow TM$ be a dilation vector field. A vector field \mbox{$X: M \rightarrow TM$} is called homogeneous of degree \mbox{$\xi\in\R$} with respect to the vector field $D$ if it satisfies the following relation:
\begin{equation}
\label{Homog_gen}
[D,X] = \xi X.
\end{equation}
\end{definition}
To derive scaling laws for the self-triggered times, we consider the standard dilation vector field~(\ref{gen_dilation_vf}). The trajectory of the vector field $D$ starting at the initial condition $x_0\in\R^n$ corresponds to the standard open ray $\{\lambda x_0: \lambda\in\R^+\}$. Using the inherent symmetries of homogeneous vector fields, a scaling law was developed in~\cite{tabuada10_tac} under the following triggering condition:
\eq{
\label{old_trig_cond}
\Gamma(\z(\tau,z_0)) = \z^T(\tau,z_0) M \z(\tau,z_0) =0, \qquad M = \left[\begin{array} {cc} -\sigma^2 I & 0 \\ 0 & I\end{array}\right], 
}
where $\sigma\in\R^+$ represents a design parameter trading control performance for number of executions. 
\begin{theorem}\cite{tabuada10_tac}
\label{ScaleHom_gen}
Let $\dot{\mathrm{z}} = Z(\mathrm{z})$ be a dynamical system homogeneous of degree $\xi\in\R$ with respect to the standard dilation vector field. The inter-execution time \mbox{$\tau:\R^{2n} \rightarrow \R^+\cup\{\infty\}$} implicitly defined by the execution rule~\eqref{old_trig_cond} 
scales according to:
\begin{equation}
\label{scal_eq_constant}
\tau(\lambda z_0) = \lambda^{-\xi} \tau(z_0), \qquad \lambda>0
\end{equation}
where $z_0\in\R^{2n}$ represents any point in the extended state space. 
\end{theorem}
In the remainder of this section we extend these results for non-homogeneous systems and non-homogeneous triggering conditions. 
\subsection{Homogenization of dynamical systems}
\label{homog_subsec}
We start describing how smooth nonlinear systems can be rendered homogeneous by introducing an auxiliary state variable. The procedure is well known for polynomial systems (see for example~\cite{baillieul80}). For instance, the following system:
\eqs{
\label{orig_sys_example}
\dot{\mathrm{x}}_1& = &\mathrm{x}_1 \mathrm{x}_2 + \mathrm{x}_2\notag\\
\dot{\mathrm{x}}_2& = &\mathrm{x}_1
}
can be rendered homogeneous of degree $1$ with respect to the standard dilation by adding an auxiliary variable $w\in\R$ satisfying $\dot{\w} = 0$:
\eqs{
\label{aux_sys_example}
\dot{\mathrm{x}}_1& = & \mathrm{x}_1\mathrm{x}_2+ \mathrm{x}_2 \mathrm{w} \notag\\
\dot{\mathrm{x}}_2& = & \mathrm{x}_1 \mathrm{w}  \\
\dot{\mathrm{w}}  & = & 0. \notag
}
Trajectories for the original system~(\ref{orig_sys_example}) can be recovered from the trajectories of the auxiliary system~(\ref{aux_sys_example}) provided that $\mathrm{w}(0)=1$. Formally speaking, we say that the vector field
\eq{X = \left[\begin{array}{ccc} x_1 x_2+x_2\\x_1\end{array}\right]
}
is $\varphi$-related with the vector field
\eq{S = \left[\begin{array}{ccc} x_1 x_2+x_2 w\\x_1 w\\0\end{array}\right]
}
with $\varphi(x)=(x,1)$, since $T\varphi \cdot X = S \circ \varphi$. This idea was exploited in our previous work~\cite{tabuada10_tac} to derive scaling laws for polynomial systems using Theorem~\ref{ScaleHom_gen}. We generalize this homogenization procedure to smooth nonlinear systems in the following lemma\footnote{Although we believe that such result might have appeared in the literature, we were not able to find it. For completeness, we include a proof.}. 
\begin{lemma}
\label{homogenization_lemma}
For any vector field
$$X(x)=\sum_{i=1}^n f_i(x) \frac{\partial}{\partial x_i},$$
there exists a vector field $S$ $\varphi$-related to $X$, with $\varphi(x)=(x,1)$, homogeneous of degree $\xi\in\R$ with respect to the standard dilation vector field, and given by:
\eq{S(x,w) = \sum_{i=1}^n w^{\xi+1} f_i (w^{-1}x) \frac{\partial}{\partial x_i}.}
\end{lemma}
\begin{IEEEproof}
We first prove that $S$ is homogeneous of degree $\xi$. The standard dilation vector field for $S$ is given by:
\eq{
\label{gen_dilation_vf_ext}
D = \sum_{i=1}^n x_i \frac{\partial}{\partial x_i}+w \frac{\partial}{\partial w}.
}
We compute the Lie bracket between $D$ and $S$ to verify condition~(\ref{Homog_gen}):
\eqs{
[D,S] &=& \frac{\partial S}{\partial x} D -  \frac{\partial D}{\partial x} S  = 
\left[\begin{array}{cc} \frac{\partial}{\partial x} \left(w^{\xi+1} f\left(w^{-1}x\right)\right) &  \frac{\partial}{\partial w} \left(w^{\xi+1} f\left(w^{-1}x\right)\right)\\ 0 & 0\end{array}\right] \left[\begin{array}{cc} x \\  w\end{array}\right]-S \notag\\
&=&\left[\begin{array}{cc}  w^{\xi+1} \frac{\partial f\left(x\right)}{\partial x} \frac{\partial (w^{-1}x)}{\partial x} x + (\xi+1) w^\xi f(w^{-1}x) w +w^{\xi+1} \frac{\partial f(x)}{\partial x}
\frac{\partial (w^{-1}x)}{\partial w} w 
\\  0\end{array}\right] - S \notag\\
&=&\left[\begin{array}{cc}   (\xi+1) w^{\xi+1}  f(w^{-1}x)   
\\  0\end{array}\right] - S = \xi S
}
since $\frac{\partial (w^{-1}x)}{\partial x} x =-\frac{\partial (w^{-1}x)}{\partial w} w $.
Moreover, the vector fields $X$ are $S$ are $\varphi$-related since~(\ref{phi_rel_condition}) is satisfied:
\eq{
T\varphi \cdot X = \left[\begin{array}{cc} f(x)\\ 0\end{array}\right] = S (x,1) = S \circ \varphi. 
}
\end{IEEEproof}
It is easy to see that the homogenization procedure for polynomial systems represents a particular case of Lemma~\ref{homogenization_lemma}. Notice that the auxiliary system is not defined for $w=0$. Nonetheless, it will be shown in section~\ref{iso_mnf_sect} that the $2n+1$ dimensional model is only used for points $(z,w)$ with $w\neq 0$. By embedding the original $2n$-dimensional dynamical system into a $2n+1$ dimensional homogeneous system  we can derive similar scaling laws for the inter-execution times, as explained in the next subsection. 


%

\subsection{Homogenization of triggering conditions}
\label{sec_homog_trig}
We also extend the class of triggering conditions to be considered. In our previous work~\cite{tabuada10_tac} we introduced the simple event-triggering condition~\eqref{old_trig_cond}.
However, in general one might require a different triggering condition. For instance, another possible condition is $\dot{V} = 0$, in order to guarantee $\dot{V}\leq 0$ for all $t\in\R^+$.

Theorem~\ref{ScaleHom_gen} leverages the fact that the triggering condition~\eqref{old_trig_cond} is quadratic in $z$. Consider now a triggering condition of the form $\Gamma(\z(t,z))=0$, where $\Gamma$ is a homogeneous function of degree $\vartheta$, and the trajectory $\z$ is governed by a vector field homogeneous of degree $\xi$. Then, the scaling laws for the inter-execution times are the same regardless of the degree of homogeneity of the function $\Gamma$. Indeed, note that the triggering condition $\Gamma$ satisfies:
\eq{
\Gamma(\z(t,\lambda z_0)) = \Gamma(\lambda \z(\lambda^{\xi} t,z_0)) = \lambda^{\vartheta} \Gamma(\z(\lambda^{\xi}t,z_0))
}
where the first equality is a property of homogeneous flows (see~\cite{kawski95} for details). Since $\lambda^{\vartheta}\neq 0$, whenever
\eq{
\label{exp_cond}
\lambda^{\vartheta} \Gamma(\z(\lambda^{\xi}t,z_0)) = 0
}
holds, the triggering condition
\eq{
\label{orig_cond}
\Gamma(\z(t,\lambda z_0)) = 0
}
is fulfilled, and vice versa. In other words, the time instants $\lambda^{\xi} t$ at which~\eqref{exp_cond} is satisfied coincide with the time instants $t$ at which~\eqref{orig_cond} is satisfied.
Hence, the same scaling laws hold whenever the triggering condition is a homogeneous function. 
  
Moreover, a similar procedure as the one used for vector fields can be exploited to render a more general class of triggering conditions homogeneous. Starting with a nonhomogeneous triggering condition $\Gamma(z)=0$ we can construct a homogeneous triggering condition of degree $\vartheta$ by means of an auxiliary variable $w$:
\eq{
w^{\vartheta+1}\Gamma(w^{-1}z) = 0
}
that is equivalent to the original triggering condition $\Gamma(z)=0$ for $w=1$. By working in a $2n+1$ dimensional space we are able to construct symmetries for the inter-execution times. These two homogenization procedures for the triggering condition and the vector field are summarized in the following theorem.
\begin{theorem}
\label{super_mega_thm}
Let \mbox{$\tau:\R^{2n} \rightarrow \R^+\cup\{\infty\}$} be the inter-execution times for the system $\dot{\z} = Z(\z)$
implicitly defined by $\Gamma(\z(\tau(z),z))=0$ for a point $z \in \R^{2n}$, and let  
\mbox{$\tilde{\tau}:\R^{2n+1} \rightarrow \R^+\cup\{\infty\}$} be the inter-execution times for the system:
\eq{
\label{homogenized_system}
\left[\begin{array}{c} \dot{\z} \\ \dot{\w} \end{array}\right] = S(\z,\w) = \left[\begin{array}{c} \w^{\xi+1} Z(\w^{-1}\z) \\ 0 \end{array}\right]
}
implicitly defined by $\tilde{\Gamma}(\z,\w):=\Gamma(\w^{-1}\z(\tilde{\tau}(z,w),z,w))=0$, for a point $(z,w)\in\R^{2n+1}$. Then, $\tau$ and $\tilde{\tau}$ are related according to:
\eq{
\label{final_thm_hom}
\tau(z) = \lambda^\xi  \tilde{\tau}(\lambda z, \lambda) \qquad \forall \lambda>0.}
\end{theorem}
\begin{IEEEproof}
We first note that $Z:M\rightarrow TM$ is $\varphi$-related to $S:N\rightarrow TN$ with $\varphi(z)=(z,1)$. 
Thus the following diagram commutes:
\begin{equation}
\label{comm_diag_app}
\begin{diagram}
\R  & \lTo^{\Gamma}         & M & \rTo^Z         & TM\\
    & \luTo_{\tilde{\Gamma}}        & \dTo^{\varphi}     &   & \dTo_{T\varphi}\\  
    &  & N & \rTo^S & TN
\end{diagram}
\end{equation}
since the following two conditions are satisfied:
\begin{enumerate}
\item $T\varphi \cdot Z = S \circ \varphi$.
\item $\tilde{\Gamma} \circ \varphi = \Gamma$.
\end{enumerate}
That is, in order to study the evolution of the map $\Gamma(z)$ under the vector field $Z$ we can alternatively study the evolution of the map $\Gamma(w^{-1}z)$ under the dynamics of $S$.
Let $\s$ be the flow for the vector field $S$. Since $Z$ and $S$ are $\varphi$-related, the corresponding flows satisfy~\cite{lee}:
\eq{\varphi \circ \z(t,z) = \s (t,\varphi(z)), \quad z\in M
}
and from commutativity of Diagram~\ref{comm_diag_app} we can conclude:
\eq{\Gamma \circ \z(t,z) = \tilde{\Gamma} \circ \varphi \circ  \z(t,z) = \tilde{\Gamma} \circ \s(t,\varphi(z)).
}
Since the evolution of both triggering conditions $\Gamma \circ \z$ and $\tilde{\Gamma} \circ \s \circ \varphi$ is identical, they generate the same inter-execution time for any point $z\in M$:
\eq{
\tau(z) = \tilde{\tau}(z,1).
}
Finally, since the vector field $S$ is homogeneous of degree $\xi$, times scale according to Equation~\eqref{scal_eq_constant}:
\eq{\tau(z) = \tilde{\tau}(z,1) = \lambda^\xi \tilde{\tau}(\lambda z,\lambda) \qquad \forall \lambda>0.
}
\end{IEEEproof}
Notice that this result does not rely on the fact that $\dot{\e} =-\dot{\x}$, but rather on the fact that the augmented vector field $S$ is homogeneous. In the following example we show how the previous theorem can be applied to study the evolution of the trajectories of a nonlinear system.
\begin{example}
Consider a simple unforced pendulum described by:
\eqs{
\label{orig_pend}
\dot{\x}_1 &= & \x_2\notag\\
\dot{\x}_2 &= & -\frac{g}{l}\sin\left(\x_1\right)-\frac{k}{m}\x_2
}
where $x_1$ corresponds to the angle of the pendulum (with respect to a vertical line), $x_2$ its angular velocity, $m$ is the mass of the pendulum, $l$ the length of the rod, $k$ is the coefficient of friction and $g$ the gravitational constant. We are interested in analyzing the time at which the angle of the pendulum reaches $\pi/6$, starting with initial condition $(x_{1_0},x_{2_0})$. Hence, we are looking for the time $\tau$ at which the following equality is satisfied:
\eq{
\label{orig_trig_example}
\Gamma(\x_1(\tau,x_{1_0},x_{2_0})) := \x_1(\tau,x_{1_0},x_{2_0}) - \pi/6 = 0.}
We homogenize the system as described before (for instance, with degree of homogeneity $\xi=1$):
\eqs{
\label{homog_pend}
\dot{\x}_1 &= & \x_2 \w\notag\\
\dot{\x}_2 &= & -\frac{g}{l}\sin\left(\w^{-1}\x_1\right)\w^2-\frac{k}{m}\x_2 \w\\
\dot{\w} & = & 0.\notag
}
Likewise, equality~\eqref{orig_trig_example} becomes:
\eq{
\tilde{\Gamma}(\x_1(\tilde{\tau},x_{1_0},x_{2_0},w_0)) := \x_1(\tilde{\tau},x_{1_0},x_{2_0},w_0) - \frac{\pi}{6}\w(\tilde{\tau},x_{1_0},x_{2_0},w_0) =0.
}
The time $\tilde{\tau}$ implicitly defined by $\tilde{\Gamma}(\x_1(\tilde{\tau},x_{1_0},x_{2_0}))=0$ and the time $\tau$ implicitly defined by $\Gamma(\x_1(\tau,x_{1_0},x_{2_0}))=0$ are related according to:
\eq{
\label{relation_example}
\tau(x_{1_0},x_{2_0}) = \lambda \tilde{\tau} (\lambda x_{1_0},\lambda x_{2_0}, \lambda) \qquad \lambda>0.
}
\end{example}
As shown in this example, the applicability of the previous framework goes beyond the self-triggered control problem and can be used to describe other temporal aspects of the trajectories of nonlinear systems. In the next section we describe how equation~\eqref{relation_example} can be used to determine the time at which equality~\eqref{orig_trig_example} is satisfied.

\section{Isochronous manifolds}
\label{iso_mnf_sect}
In this section we exploit the results of the previous section to derive a self-triggering strategy that emulates an event-triggered implementation. In Section~\ref{sec_scal_laws} scaling laws were derived for the inter-execution times of dynamical systems. Such scaling laws describe the evolution of the times implicitly defined by a triggering condition. Once the inter-execution time for a point is known, times for any point lying on the same ray can be derived by means of the scaling laws. Hence, in order to compute the inter-execution times for the whole extended state space it is necessary to obtain the inter-execution times for at least one point in each ray. For instance, it would be enough to find a lower bound $\tau^*$ for the times on a sphere (since every ray intersects a sphere), and then extend the times along rays. This was the approach used in our previous work~\cite{tabuada10_tac}. Different methods are available to compute a lower bound $\tau^*$ for the inter-execution times on a sphere (see for instance~\cite{tabuada07}). Following this procedure, inter-execution times for the self-triggering strategy can be computed using~(\ref{scal_eq_constant}):
\eq{
\label{old_self}
\tau(\z(t_i)) = \lambda^{-\xi} \tau^*, \qquad \quad\text{with} \quad \lambda = \frac{|z(t_i)|}{r} 
}
where we have identified the point $\lambda z_0$ in~(\ref{scal_eq_constant}) with the last measurement of the extended state $\z(t_i) = (\x(t_i),0)$, the left hand side of~(\ref{scal_eq_constant}) with the time for $\z(t_i)$, $z_0$ with a point in a sphere of radius $r$, and $\tau(z_0)$ with $\tau^*$, the lower bound for such sphere. Unfortunately, such procedure introduces two sources of conservativeness: 
\begin{itemize}
	\item The existing methods to compute a stabilizing period for a given set (also known in the literature as maximum allowable transfer interval, MATI) are rather conservative~\cite{nesic2009explicit,tabuada07}, since in general they are based on the Lipschitz constant of the vector field.
	\item The stabilizing period represents a lower bound for the inter-execution times of all the points in the sphere. Since times might vary drastically along the sphere (or along any given set), this lower bound will be overly conservative for some points in the sphere. In other words, the evolution of the times across rays is neglected in this procedure.
\end{itemize}
	
In the next subsections we overcome these two drawbacks through the use of isochronous manifolds.

\subsection{Existence of isochronous manifolds}
\label{sec_iso_mnf}
As mentioned before, the sphere was selected in~\cite{tabuada10_tac} as the reference manifold (i.e., the set where a lower bound for the inter-execution times is computed) for simplicity. However, there is no reason to expect that such times are constant for the points in the sphere. Ideally, we would like to identify manifolds where the inter-execution times remain constant. Such manifolds will be termed as isochronous manifolds. If these manifolds can be computed in an exact manner, the inter-execution times generated by the self-triggering strategy will coincide with the event-triggered times.

It is important to notice that we are interested in isochronous manifolds that allow us to derive self-triggering conditions for the whole operating region through the scaling laws of Section~\ref{sec_scal_laws}. In other words, in order to be of any use, an isochronous manifold should intersect \emph{every} ray at least at one point. We prove the existence of this class of isochronous manifolds for homogeneous systems in the following proposition.
\begin{proposition}
\label{existence_iso_mnf}
Consider the  dynamical system $\dot{\z} = Z(\z)$ homogeneous of degree $\xi>0$, defined on a $2n$-dimensional manifold. Let $\tau:\R^{2n} \rightarrow \R^+ \cup \{\infty\}$ be implicitly defined by $\Gamma(\z(\tau(z_0),z_0))=0$, where $\Gamma$ is a homogeneous function. If the following two assumptions are satisfied for any $z_0$ in a sphere of radius $r$ (for some $r\in \R^+$):
\eq{\label{assump1_prop} \Gamma(\z(0,z_0))\neq 0,}
\eq{\label{assump2_prop} \exists t \in \R^+: \Gamma(\z(t,z_0))= 0}
then, for any given time $t_*\in \R^+$ the set $\Omega_{t_*}$ defined by:
\begin{equation}
\label{eq_iso_mnf}
\Omega_{t_*} = \{z_*\in \R^{2n}: \tau(z_*) = t_*\}
\end{equation}
is a manifold of dimension $2n-1$.
\end{proposition}
\begin{IEEEproof}
As explained in Section~\ref{sec_homog_trig}, the inter-execution time $\tau$ defined by the triggering condition scales according to:
\begin{equation}
\label{scal_const_degree}
\hspace{1.8cm}\tau(\lambda z_0) = \lambda^{-\xi} \tau(z_0), \hspace{0.4cm}  \forall \lambda>0,\xi >0. \\
\end{equation}
According to~(\ref{scal_const_degree}), inter-execution times vary from $0$ to $\infty$ as $\lambda$ varies from $\infty$ to $0$. Hence, for any given time $t_* \in \R^+$, there exists a point $z_*$ in each ray such that $\tau(z_*) = t_*$. Moreover, equation~(\ref{scal_const_degree}) implies that there does not exist two different points lying on the same ray with the same inter-execution time. The union of all those points $z_*$ constitutes an isochronous manifold for the homogeneous system. 

We briefly sketch now the proof that $\Omega_{t^*}$ is indeed a manifold. We use $\mathrm{d}$ to denote the exterior derivative of a function. We first observe that $\Omega_{t_*}$ can be equivalently defined as the zero level set of the map $\Gamma(\mathrm{z}(t_*,z))$ for the chosen $t_*$. 
Since $\Gamma$ is homogeneous we have:
\eq{\Gamma(\lambda z) = \lambda^{\vartheta} \Gamma(z), \quad \forall \lambda>0}
where $\vartheta$ is the degree of homogeneity of $\Gamma$. Then,~\eqref{assump1_prop} implies that the triggering condition never vanishes for any $z \neq 0$:
\eq{
\Gamma(z) = \Gamma\left(\frac{|z|}{r}\frac{rz}{|z|} \right) = \left(\frac{|z|}{r}\right)^{\vartheta}\Gamma\left(\frac{rz}{|z|} \right)\neq 0
}
as the point $\frac{rz}{|z|}$ belongs to a sphere of radius $r$. Since $\Gamma$ is a homogeneous function
it satisfies:
\eq{
\label{prop_partials}
{\cal L}_D \Gamma = (\xi +1) \Gamma
}
with $D$ being the standard dilation vector field, as defined in~\eqref{gen_dilation_vf}. Since the Lie derivative can also be expressed as:
\eq{
\label{ext_der}
\mathrm{d}\Gamma (D) = {\cal L}_D \Gamma. 
}
and $D$ only vanishes at the origin, we conclude from~\eqref{prop_partials} and~\eqref{ext_der} that $\mathrm{d}\Gamma$ only vanishes when $\Gamma=0$, \textit{i.e.}, at $z=0$; therefore, $\Gamma$ is a submersion for every $z\ne 0$. Since $\mathrm{z}(t_*,z)$ is a diffeomorphism for every $t_*\in\R$, the corresponding tangent map is an isomorphism for every \mbox{$z\in \R^n$}. Hence, we conclude that $\Gamma(\mathrm{z}(t_*,z))$ is a smooth submersion from which follows that its zero level set is a smooth manifold of codimension 1~\cite{lee}.
\end{IEEEproof}
Let us comment now on the two assumptions~\eqref{assump1_prop} and~\eqref{assump2_prop}. Assumption~\eqref{assump1_prop} will be satisfied by the triggering conditions that are usually of interest. Indeed, should the triggering condition $\Gamma(\z(0,z_0))$ be $0$ for some $z_0$, the inter-execution time for such $z_0$ will be 0, which is undesirable from a practical point of view.
For instance, for the triggering condition~\eqref{old_trig_cond} mentioned in Section~\ref{sec_scal_laws}, we have $\e(0,z_0)=0$, hence assumption~\eqref{assump1_prop} holds:
\eqs{
\Gamma(\z(0,z_0)) &=& \z(0,z_0)^T M \z(0,z_0) = \e(0,z_0)^T \e(0,z_0) - \sigma^2 \x(0,z_0)^T \x(0,z_0) \notag\\
&=& - \sigma^2 \x(0,z_0)^T \x(0,z_0) < 0 \quad \forall z_0\neq 0.
}
Likewise, it can be easily checked that the triggering condition in~\eqref{trig_condition} satisfies this condition. 

Assumption~\eqref{assump2_prop} guarantees that there does not exist states where the inter-execution times are $\infty$. In general, there might exist stable manifolds where no updates in the control input are needed to steer the system to the equilibrium point. In this case, the isochronous manifold would not intersect every homogeneous ray, and therefore it would not be useful to derive self-trigger formulas. For those states lying in stable manifolds, the scaling law~(\ref{scal_const_degree}) implies that times will be $\infty$ for all points lying on the rays intersecting the stable manifold. And vice versa, since all those points have $\infty$ as the next inter-execution time, any stable manifold for homogeneous systems needs to be composed of homogeneous rays. Thus, Assumption~\eqref{assump2_prop} guarantees that times on each ray assume values in $]0,\infty[$.


The isochronous manifold can be used instead of the sphere in order to compute the inter-execution times for the whole operating region, as there exists a unique intersection point between each ray and the manifold. In the next section we provide techniques for the computation of such manifolds.

\subsection{Approximating isochronous manifolds}
\label{sec_self}
The isochronous manifold $\Omega_{t_*}$ implicitly defined by equation~(\ref{eq_iso_mnf}) represents the set of all extended states $z_*$ that satisfy the triggering condition at time $t_*$.  Alternatively, the isochronous manifold can be defined as:
\begin{equation}
\label{def_iso_mnf}
\Omega_{t_*} = \left\{z_* \in \R^n \big \vert \Gamma(\z(t_*,z_*)) =0\right\}.
\end{equation}
The explicit computation of $\Omega_{t_*}$ would require the knowledge of the flow, which is in general unknown. Herein we develop a technique to approximate the isochronous manifold. 
Since we  are only interested in the evolution of the triggering condition $\Gamma (\mathrm{z}(t,z_0))=0$, it would be sufficient to construct a differential equation for $\Gamma$.  For ease of exposition and without loss of generality, at this point we assume that $\Gamma (\z(0,z_0)) \leq 0$ for any $z_0$ in the operating region (which is the case for the triggering conditions~\eqref{old_trig_cond} and $\dot{V} =0$).
In~\cite{tabuada07} a first order nonlinear differential inequality of the form \mbox{$\frac{d}{dt}\Gamma \leq \nu(\Gamma)$} was derived, whose solution just depended on the norm of the initial extended state $z_0$. Instead, we derive a $p$-th order differential equation, whose solution collects more information from the extended state. Moreover, in order to obtain a closed-form expression for its solution, we focus on linear differential equations describing  the evolution of \mbox{$\Gamma$}. Clearly, $\Gamma$ will be described by a linear differential equation only under special circumstances, see \eg~\cite{levine1986nsi}. However, the evolution of \mbox{$\Gamma$} can be \emph{bounded} by a linear differential equation. This amounts to constructing a set of coefficients $\chi_i\in\R$ satisfying:
\begin{equation}
\label{ineq_lies}
{\mathcal L}_Z^p \Gamma(z) \leq \sum_{i=0}^{p-1} \chi_i {\mathcal L}_Z^i \Gamma(z)
\end{equation}
with $Z$ representing the closed-loop vector field and for any extended state $z$ in an invariant region of interest $M\subseteq \R^{2n}$, as stated in the next result.
\begin{lemma}
\label{lemma_approx_trig}
Consider a vector field $Z:M\to TM$ and a map $\Gamma:M\to \R$. For every set of coefficients $\chi_0,\chi_1,\hdots\chi_{p-1}\in \R$ satisfying:
\begin{equation}
\label{ineq_lies_lemma}
{\mathcal L}_Z^p \Gamma(z) \leq \sum_{i=0}^{p-1} \chi_i {\mathcal L}_Z^i \Gamma(z) \qquad \forall z\in M,
\end{equation}
the following inequality holds for all $z_0\in M$ and for all $t\in \R_0^+$ for which the solution $\z$ of $Z$ is defined:
\begin{equation}
\label{comp_traj}
\Gamma(\z(t,z_0))\le \mathrm{y}_1(t,y_0)
\end{equation}
where $\mathrm{y}_1$ is the first component of the solution of the linear differential equation:
\begin{eqnarray}
\dot{\mathrm{y}}_i &=&\mathrm{y}_{i+1}, \qquad i=1,\ldots,p-1\\
\label{def_yp}
\dot{\mathrm{y}}_p &=&\sum_{i=0}^{p-1} \chi_i \mathrm{y}_i
\end{eqnarray}
with initial condition:
\eq{
\label{eq_ini_conditions}
y_{0}=(\Gamma(z_0), {\mathcal L}_Z \Gamma(z_0), \hdots, {\mathcal L}_Z^{p-1} \Gamma(z_0))^T.
}
\end{lemma}
\begin{IEEEproof}
The differential inequality~\eqref{ineq_lies_lemma} can be rewritten as:
\eqs{
\label{def_bp}
\dot{\mathrm{b}}_i & = &\mathrm{b}_{i+1} \qquad i=1,\ldots,p-1\\
\dot{\mathrm{b}}_p &\leq &\sum_{i=0}^{p-1} \chi_i \mathrm{b}_i
}
with initial condition $b_{0}=(\Gamma(z_0), {\mathcal L}_Z \Gamma(z_0), \hdots, {\mathcal L}_Z^{p-1} \Gamma(z_0))^T$.
By the Comparison Lemma (see for instance~\cite{khalil2002ns}) we have $\mathrm{b}_p(t,b_0)\le \y_p(t,y_0)$ since $b_0=y_0$ and $\dot{\mathrm{b}}_p(t,b_0)={\mathcal L}_Z^p \Gamma(\z(t,z_0)) \le \sum_i \chi_i {\mathcal L}_Z^i \Gamma (\z(t,z_0))=\dot{\y}_p(t,y_0)$, in virtue of~(\ref{ineq_lies_lemma}). 
We apply the Comparison Lemma successively to conclude that each $\y_i$ bounds each $\mathrm{b}_i$:
\begin{equation}
\mathrm{b}_i(t,b_0) \leq \mathrm{y}_i(t,y_0), \qquad i=1,\hdots,p, \quad t\geq 0.
\end{equation}
Hence, $\y_1$ represents a bound for the evolution of $\Gamma$:
\eq{
\Gamma (\z(t,z_0)) = \mathrm{b}_1(t,b_0) \le \mathrm{y}_1(t,y_0)\notag
}
\end{IEEEproof}
For convenience, at this point we define the following map:
\eq{
\mu^p (z):= (\Gamma(z), {\mathcal L}_Z \Gamma(z), \hdots, {\mathcal L}_Z^{p-1} \Gamma(z))^T.
}
Using Lemma~\ref{lemma_approx_trig} we can see that the minimum time satisfying $\mathrm{y}_1 (t,y_0)=0$ lower bounds the minimum time satisfying $\Gamma(\z(t,z_0))=0$, since $\Gamma(\mathrm{z}(0,z_0))\leq 0$. Thus, the inter-execution time $\tau^{\downarrow}$ defined by $\mathrm{y}_1(\tau^{\downarrow},\mu^p(z))=0$ represent a lower bound for the inter-execution time $\tau$ defined by the triggering condition $\Gamma(\z(\tau,z_0))=0$. We can interpret this lemma as a way to construct a linear dynamical system $\dot{\mathrm{y}} = A_p \mathrm{y}$ that bounds the evolution of $\Gamma$, with $A_p$ defined as:
\begin{equation}
\label{def_Ap}
A_p = \left[\begin{array}{cccccc} 
0 &1& 0& \ldots &0& 0\\
0& 0& 1& \ldots& 0& 0 \\ 
\vdots& & & \ddots &&\vdots \\
0 & 0 & 0 &\hdots & 1 & 0\\
0 & 0 & 0 &\hdots & 0 & 1\\
\chi_0 & \chi_1 & \chi_2 &\hdots &\chi_{p-2} & \chi_{p-1} \end{array} \right]. 
\end{equation}
Using the trajectories of this linear system we rewrite the triggering condition $\mathrm{y}_1(\tau^{\downarrow},\mu^p(z))=0$ as:
\begin{equation}
\label{new_trig_cdn}
\mathrm{y}_1 (\tau^{\downarrow},\mu^p(z)) = \left[\begin{array} {cccc}1 & 0 & \hdots & 0\end{array}\right] \exp^{A_p \tau^{\downarrow}} \left[\begin{array} {c} \Gamma(z) \\ 
({\mathcal L}_Z\Gamma)(z) \\ \vdots \\({\mathcal L}_Z^{p-1} \Gamma)(z) \end{array}\right]=0.
\end{equation}
It is clear that by fixing $\tau^{\downarrow}=t_*$ in equation~(\ref{new_trig_cdn}) we obtain an equation describing the set of states whose inter-execution times are lower bounded by $t_*$. Hence the set $\Omega^{\downarrow}_{t_*}$ defined by:
\begin{equation}
\label{aprox_iso_mnf}
\Omega^{\downarrow}_{t_*} = \{z_*\in\R^{2n}: \mathrm{y}_1 (t_*,\mu^p(z_*))=0\}
\end{equation}
represents a bound for the isochronous manifold $\Omega_{t_*}$ in the following sense:
$$z_* \in \Omega^{\downarrow}_{t_*} \quad \implies \quad \tau(z_*)\geq t_*.$$
As mentioned before, we will use this bound to compute the inter-execution times through the scaling laws derived in section~\ref{sec_scal_laws}. 

\subsection{A self-triggering formula}
\label{sec_self_trig_formula}
In order to apply the scaling law~(\ref{scal_eq_constant}) with $\Omega^{\downarrow}$, we first find the intersection between the homogeneous rays and $\Omega^{\downarrow}$. Since we are searching for a self-triggering condition to be applied online, it is desirable to have a closed-form expression for those intersecting points. Towards this objective we state a simple lemma describing a useful property of the Lie derivative of a homogeneous function along homogeneous vector fields.
\begin{lemma}
\label{homog_lemma}
Consider a map $\Gamma:M \rightarrow \R$ homogeneous of degree $\vartheta$ and a vector field \mbox{$Z:M \rightarrow TM$} homogeneous of degree $\xi\in\R$. Then, the $k$-th Lie derivative of $\Gamma$ along $Z$ is homogeneous of degree $\vartheta+k \xi$:
\begin{equation}
{\mathcal L}_Z^k \Gamma(\lambda z) = \lambda^{\vartheta+k \xi} {\mathcal L}_Z^k \Gamma(z), \qquad z\in M, \lambda>0.
\end{equation}
\end{lemma}
The result can be easily proven by induction. To make use of the scaling law~(\ref{scal_eq_constant}), we identify $\lambda z_0$ with a point $z_*$ in $\Omega^{\downarrow}$ and $z_0$ with $\z(t_i)$. To compute the intersection point between a ray passing through $\z(t_i)$ and $\Omega^{\downarrow}$ we substitute the equation of the ray $z_* = \lambda \z(t_i)$ in the triggering condition~(\ref{new_trig_cdn}):
\begin{eqnarray}
\label{eq_lambdas}
\hspace{0cm}\mathrm{y}_1 (t_*,\mu^p(z_*)) &\!\!\!\!=\!\!\!\!& \mathrm{y}_1 (t_*,\mu^p(\lambda\z(t_i))) 
= \left[\begin{array} {cccc}\!\!\!1 \!\!\!& 0 \!\!\!& \!\hdots \!&\!\!\! 0\!\!\!\end{array}\right] \exp^{A_p t_*} \!\!\left[\begin{array} {c} \Gamma(\lambda \z(t_i)) \\ {\mathcal L}_Z\Gamma(\lambda \z(t_i)) \\ \vdots \\{\mathcal L}_Z^{p-1} \Gamma(\lambda \z(t_i)) \end{array}\right] \notag\\
&\!\!\!\!=\!\!\!\!&\sum_{i=0}^{p-1} \beta_i(t_*,\z(t_i)) \lambda^{\vartheta+i\xi} = 0
\end{eqnarray}
with:
\begin{equation}
\label{def_beta}
\beta_i(t_*,\z(t_i)) = {\exp^{A_p t_*}}_{1(i+1)} \:\:\:{\mathcal L}_Z^i \Gamma(\z(t_i)), \qquad i=0,\hdots,p-1. 
\end{equation}
Since we are interested in the positive real solution, we can divide equation~(\ref{eq_lambdas}) by $\lambda^{\vartheta}$ and define the variable $q=\lambda^{\xi}$ to obtain the equivalent equation:
\begin{equation}
\label{final_pol_eq}
\sum_{i=0}^{p-1} \beta_i(t_*,\z(t_i)) q^i = 0.
\end{equation}
The positive real solution of this equation can be substituted into the scaling law in~(\ref{scal_eq_constant}) to obtain a self-triggering condition, as summarized in the following theorem:
\begin{theorem}
\label{final_thm}
Consider a dynamical system $\dot{\z} = Z(\z)$, homogeneous of degree $\xi>0$. Let \mbox{$\tau:\R^{2n} \rightarrow \R^+ \cup \{\infty\}$}  be the times implicitly defined by $\Gamma(\z(\tau(z),z))=0$, where $\Gamma$ is a homogeneous function.
For any set of coefficients $\chi_0,\chi_1,\hdots,\chi_{p-1}\in \R$ satisfying~(\ref{ineq_lies}), and for any $t_*\in \R_0^+$ the submanifold $\Omega^{\downarrow}_{t_*}$, defined by~(\ref{aprox_iso_mnf}), upper bounds the isochronous manifold $\Omega_{t_*}$ in the following sense:
$$z\in \Omega^{\downarrow}_{t_*}\quad\implies\quad \tau(z)\geq t_*.$$
Moreover, the function $\tau^\downarrow:M\to \R$ defined by:
\eq{
\label{final_formula}
\tau^{\downarrow}(\z(t_i)) = q t_*\\
}
with $q$ satisfying:
\eq{
\label{final_pol_eq_thm}
\sum_{i=0}^{p-1} \beta_i(t_*,\z(t_i)) q^i = 0
}
provides a lower bound for the time $\tau(\z(t_i))$, with $\beta_i(t_*,\z(t_i))$ given by~(\ref{def_beta}). 
\end{theorem}
Equation~\eqref{final_pol_eq} reflects the tradeoff between complexity and accuracy of the computation of the inter-execution times: increasing $p$ improves the approximation, but also increases the computational cost. The value $p=3$ seems to be a sensible choice since it allows us to find a simple closed-form solution for~(\ref{final_pol_eq}). Equation~\eqref{final_formula} then becomes:
\eq{
\label{2nd_ord_self}
\tau^{\downarrow} = \frac{2\beta_0}{-\beta_1+\text{sign}(\beta_0)\sqrt{\beta_1^2 - 4 \beta_2\beta_0}} \:\:t_*.
}
\begin{remark}
For non-homogeneous systems and non-homogeneous triggering conditions, a lower bound $\tau^{\downarrow}(\z(t_i))$ for the inter-execution time $\tau(\z(t_i))$ is also given by~(\ref{final_formula}), where now the coefficients $\beta_i$ are computed with the homogenized system and the homogenized triggering condition, as defined in section~\ref{sec_scal_laws}. 
\end{remark}
\begin{remark}
In the case of linear systems, the inter-execution times remain constant along homogeneous rays since the degree of homogeneity is 0, as dictated by Theorem~\ref{ScaleHom_gen}. In order to obtain self-triggering formulas through Theorem~\ref{final_thm}, the linear system has to be first rendered homogeneous of degree $\xi>0$, as explained in Section~\ref{homog_subsec}.
\end{remark}

\section{Example}
\label{sec_example}
\begin{figure}[ht]
   \begin{center}
 	\includegraphics[width=0.75\hsize]{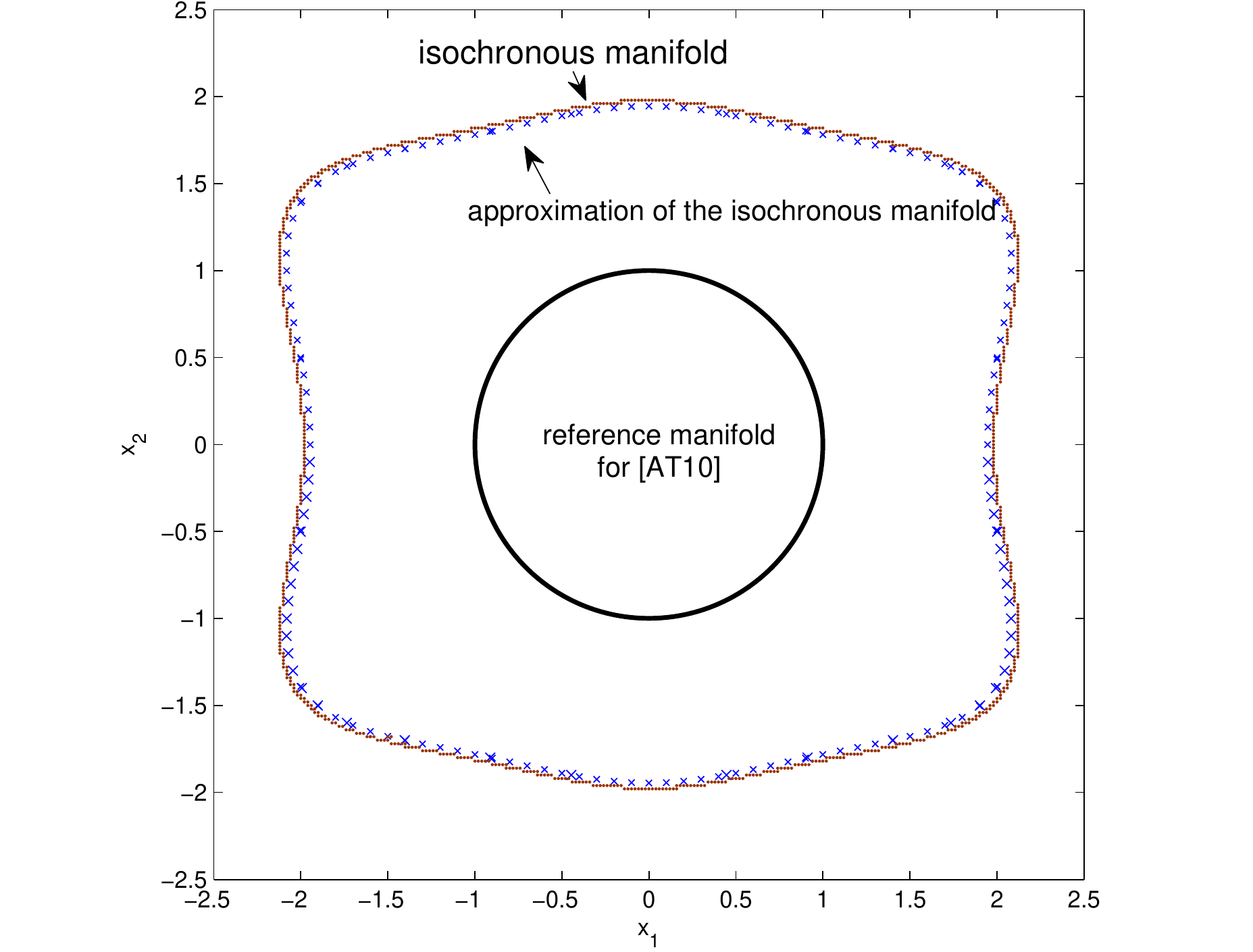}
   \end{center}
 \caption{Reference manifold for \cite{tabuada10_tac}, isochronous manifold and its approximation according to equation~(\ref{aprox_iso_mnf}).}
\label{ref_mnf}
\end{figure}

We compare the results herein developed with our previous work in~\cite{tabuada08_acc} and~\cite{tabuada10_tac}. The equations of the example in~\cite{tabuada08_acc} are:
\begin{eqnarray}
\dot{\x_1} & = & - \x_1^3 + \x_1 \x_2^2  \nonumber\\
\dot{\x_2} & = & \x_1 \x_2^2 + \u - \x_1^2 \x_2
\end{eqnarray}
with control law $u = - x_2 ^3 - x_1 x_2^2$. The operating region is a ball of radius 1 around the origin. We use the same triggering condition as in~\cite{tabuada08_acc}, given by:
\eq{
\label{trig_cond_homog_system}
|e|^2 = 0.0127^2 \sigma^2 |x|^2,\qquad \sigma\in]0,1[.
}
Using SOSTOOLS~\cite{sostools}, we obtain a set of coefficients $\chi_i$ satisfying inequality~\eqref{ineq_lies_lemma} for $p=3$: 
$$ \chi_0 = 105.970, \quad \chi_1= 0.021, \quad \chi_2 = 1.033.$$
These coefficients define the approximation of the isochronous manifold.
Figure~\ref{ref_mnf} depicts $\Omega^{\downarrow}_{t_*}$ for $t_*=1$ms, computed according to equation~(\ref{aprox_iso_mnf}), and the isochronous manifold $\Omega_{t_*}$ computed via numerical simulations. As proven in Theorem~\ref{final_thm}, the exact isochronous manifold encloses $\Omega^{\downarrow}$, since times enlarge as we approach the origin. Hence, two conclusions can be drawn from this figure. First, we notice that for the chosen value of $t_*$ the set $\Omega^{\downarrow}$ nearly coincides with $\Omega$, while there is a considerable gap between $\Omega$ and the sphere, showing that the lower bound developed in~\cite{tabuada10_tac} was not tight. Moreover, the shape of the isochronous manifold is clearly different from the sphere, implying that even if a tight bound for the times on a sphere is computed, there exist many points where the inter-execution times will be significantly conservative.

Since we are not aware of any other work addressing self-triggering strategies for nonlinear systems, in the simulations we compare the inter-execution times defined by Theorem~\ref{final_thm} against the periodic strategy, the event-triggered times generated by the triggering condition~(\ref{trig_condition}) and the self-triggering technique in~\cite{tabuada10_tac}. The period is computed as explained in~\cite{tabuada07}. 
The trajectories are similar under the 4 different implementations for all tested initial conditions.

\begin{figure}[ht]
   \begin{center}
 	\includegraphics[width=0.75\hsize]{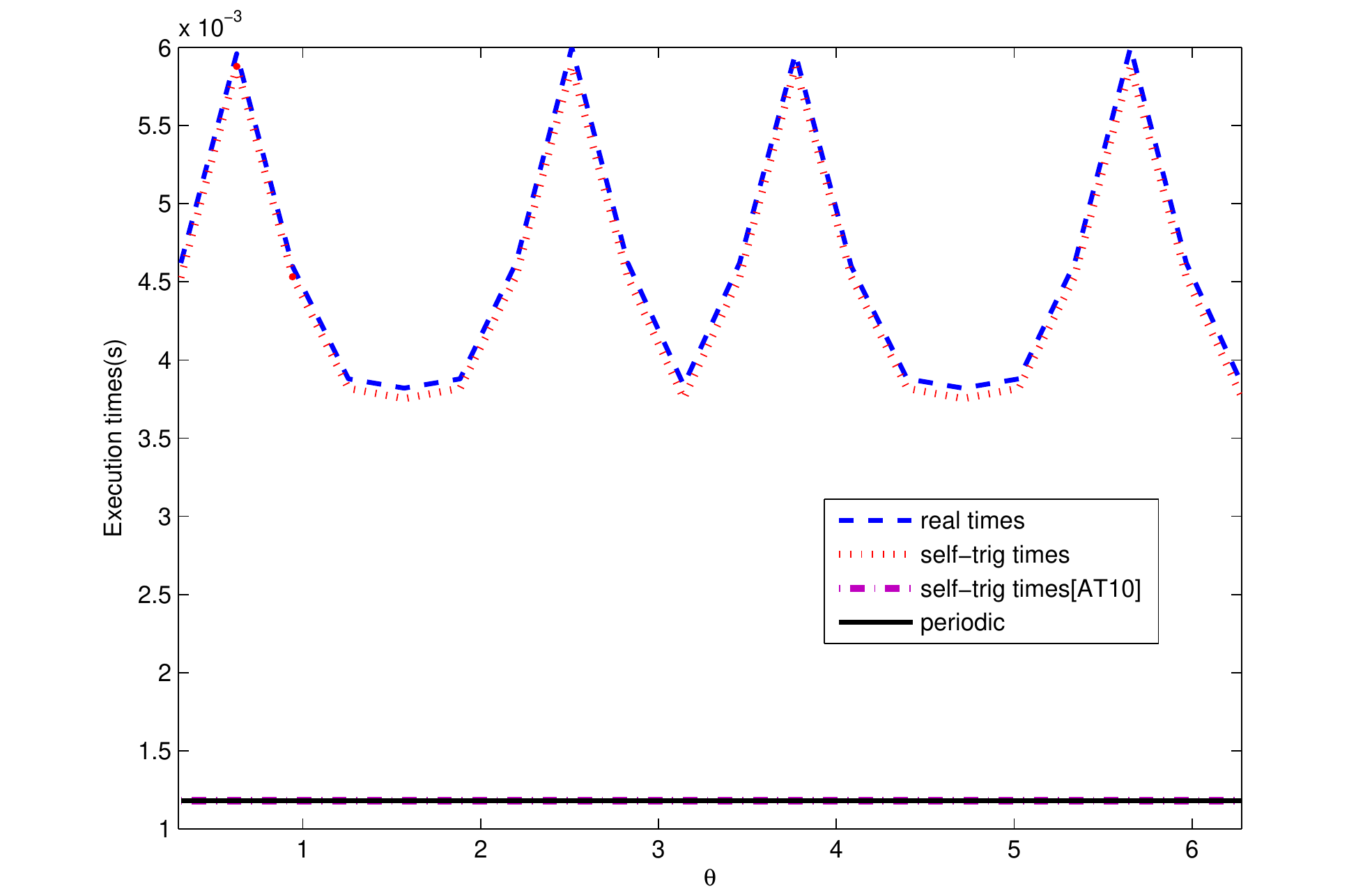}
   \end{center}
 \caption{Inter-execution times along a sphere of unitary radius.}
\label{exec_times_sphere}
\end{figure}

As our previous work in~\cite{tabuada10_tac} describes the exact evolution of the times along rays, we focus our comparison on the evolution of times across rays, which is the main topic addressed in this paper. For that purpose, we consider 20 initial conditions equally spaced along a sphere. The self-triggering condition developed in~\cite{tabuada10_tac} defines the inter-execution times just as a function of the norm of the state, hence it will produce the same inter-execution times for all points lying on the sphere. On the other hand, Theorem~\ref{final_thm} takes into account more information contained in the state. Figure~\ref{exec_times_sphere} shows the inter-execution times under the 4 strategies as a function of the position in the boundary of the sphere (that is, $x = (\cos(\theta),\sin(\theta))$, for $\theta \in [0, 2\pi[$). Table~\ref{comparison_table_acc_example} represents the average inter-execution time for the 20 points along the boundary for different values of $\sigma$ (\ie, for different degrees of performance). We can observe both in Table~\ref{comparison_table_acc_example} and Figure~\ref{exec_times_sphere} how the new proposed technique nearly matches the event-triggered times and improves significantly the times generated by the former self-triggering strategy, which generates in this case the same times as the periodic strategy. 

\floatsep 2pt plus 1pt minus 1pt
\textfloatsep=0.6\floatsep
\intextsep\floatsep
\setlength{\abovecaptionskip}{0pt}
\begin{table}
\centering
\begin{tabular}{ccccc}
\multirow{2}{*} {$\sigma$} & \multirow{2}{*}{periodic} & self-triggered & self-triggered & \multirow{2}{*}{event-triggered} \\ &&\cite{tabuada10_tac} & Thm.~(\ref{final_thm}) &\\
\hline\\
0.1 &  0.39 & 0.39 & 1.50 & 1.55\\
0.2 &  0.79 & 0.79 & 3.00 & 3.06 \\
0.3 &  1.18 & 1.18 & 4.50 & 4.58\\
\hline\\
\end{tabular}
\caption{Average time along a sphere for the~\cite{tabuada08_acc} example (in ms.)}
\label{comparison_table_acc_example}
\end{table}
\setlength{\belowcaptionskip}{0pt} 
\setlength{\abovecaptionskip}{1pt}

\floatsep 2pt plus 2pt minus 2pt
\textfloatsep=0.9\floatsep
\intextsep\floatsep
\setlength{\abovecaptionskip}{0pt}
\setlength{\belowcaptionskip}{0pt} 

\floatsep 2pt plus 2pt minus 2pt
\textfloatsep=0.9\floatsep
\intextsep\floatsep
\setlength{\abovecaptionskip}{0pt}
\begin{figure}[ht]
   \begin{center}
 	\includegraphics[width=0.75\hsize]{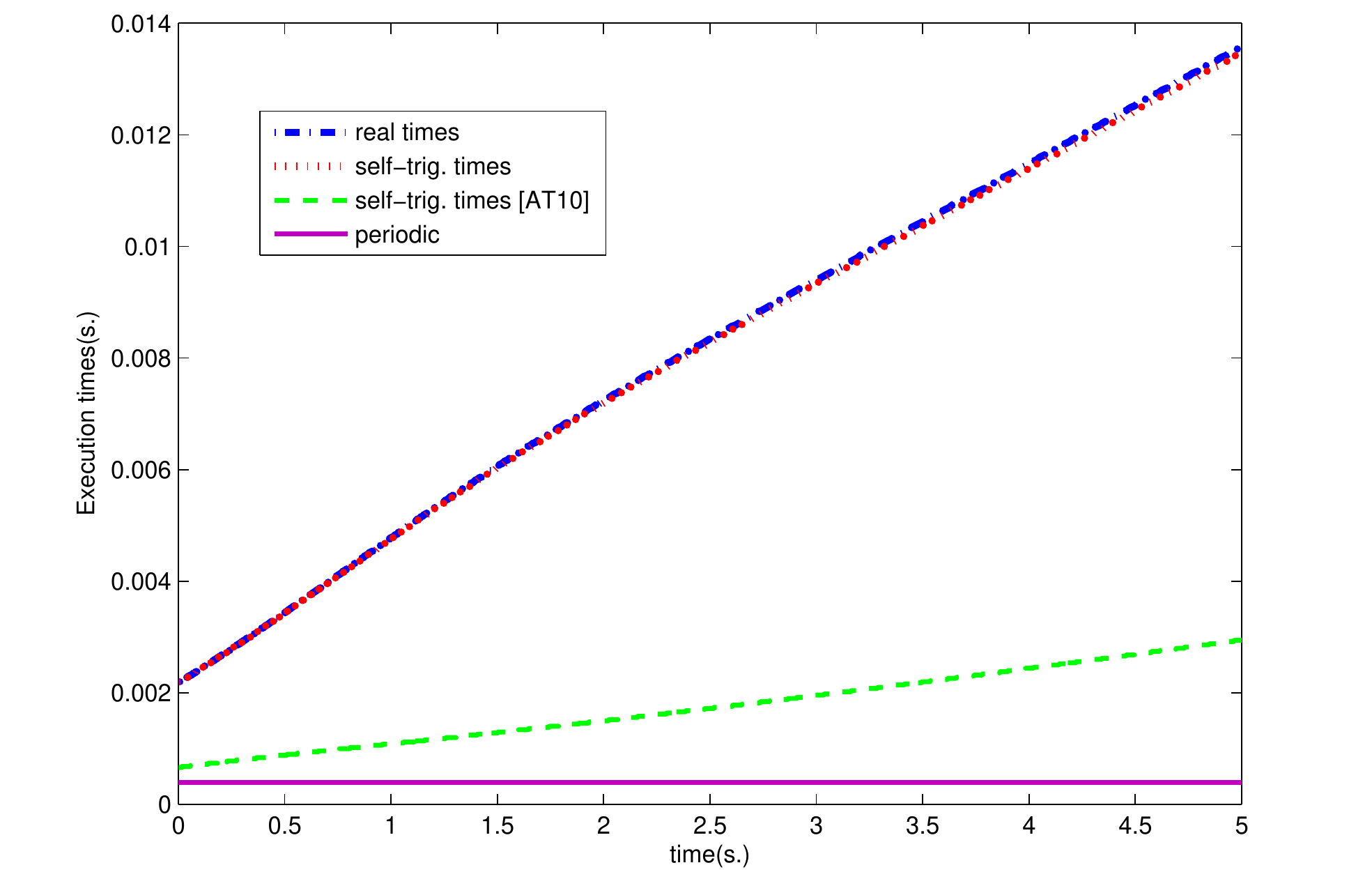}
   \end{center}
 \caption{Evolution of the inter-execution times along trajectories.}
\label{evolution_times}
\end{figure}
\setlength{\belowcaptionskip}{0pt} 
Finally, Figure~\ref{evolution_times} shows the evolution of the inter-execution times along trajectories for the 4 different implementations for a particular initial condition ($x_0=(0.4, 0.7)$) and for a simulation of 5s. Again, we can observe that the new self-triggering technique nearly tracks the evolution of the event-triggered times.

\section{Accuracy of the self-triggering technique}
\label{sec_acc}
The derived self-triggering formulas represent an accurate bound for the event-triggered times in the example in Section~\ref{sec_example}. However, there is no guarantee that 
such bounds will be accurate in general. Lemma~\ref{lemma_approx_trig} provides a constructive procedure to derive a linear model lower bounding the evolution of the triggering condition. The same procedure can be used to derive an upper bound $\tau^{\uparrow}$, by reversing the inequality~\eqref{ineq_lies}:
\eq{
\label{up_bound}
({\mathcal L}_{Z}^p \Gamma)(z,w) \geq \sum_{i=0}^{p-1} \chi_i ({\mathcal L}_{Z}^i \Gamma)(z,w).
}
Using this upper bound $\tau^{\uparrow}$, a bound for the mismatch between the self-triggered times and the event-triggered times can thus be computed:
\eq{
0 \leq \tau(\z(t_i)) - \tau^{\downarrow}(\z(t_i)) \leq \tau^{\uparrow}(\z(t_i)) - \tau^{\downarrow}(\z(t_i)).
}
While the self-triggered times generated by Theorem~\ref{final_thm} for $p\leq 3$ might be useful for some applications, they could be conservative in other cases. As mentioned before, self-triggered formulas are expected to yield larger inter-execution times
as the value of $p$ is increased, that is, when high order approximations are used in~\eqref{ineq_lies_lemma}. However, for $p>3$ no simple closed-form expression can be obtained for equation~\eqref{final_pol_eq}. To overcome this drawback, in the next subsection we develop an iterative algorithm that considers high order approximations ($p>3$) to compute a self-triggered formula.

\subsection{An iterative algorithm}

To compute the next inter-execution time $\tau(\z(t_i))$  equation~(\ref{final_pol_eq}) needs to be solved online upon the measurement of the augmented state $\z(t_i)$. Equation~\eqref{final_pol_eq} is polynomial of order $p-1$ in the variable $q$. For $p\geq4$, closed-form expressions for $\tau$ are too complicated to be evaluated online. There is a vast literature on algorithms that find approximations of the roots of polynomials. However, we seek approximate solutions that are \textit{always} lower bounding the exact root of~(\ref{final_pol_eq}), in order to obtain inter-execution times that guarantee stability of the system. For that matter, we develop an iterative algorithm that always produces lower bounds for the inter-execution times. 

For convenience we define the map:
\eq{
\delta^{\lambda}(x) := (x_1, \lambda x_2, \hdots, \lambda^{p} x_p).
}

Notice that Lemma~\ref{lemma_approx_trig} allows us to derive bounding models of different orders. We first construct two bounding models using Lemma~\ref{lemma_approx_trig}: a low order model $\underline{\y}$, with order $\underline{p}\leq 4$, and a high order model $\overline{\y}$, with order $\overline{p}\geq 4$: 
\eqs{
\label{def_ys}
\underline{\y}(t,\mu^{\underline{p}}(z)) & := & \exp^{A_{\underline{p}} t} \mu^{\underline{p}}(z)\notag\\
\overline{\y}(t,\mu^{\overline{p}}(z))  & := &  \exp^{A_{\overline{p}} t} \mu^{\overline{p}}(z)
}
The low order model will be used to compute the solutions and the high order model will serve to evaluate the accuracy of the solution, in the spirit of the Newton-Raphson method. 

\setlength{\abovecaptionskip}{0pt}
\begin{figure}[ht]
   \begin{center}
 	  \includegraphics[width=0.75\hsize]{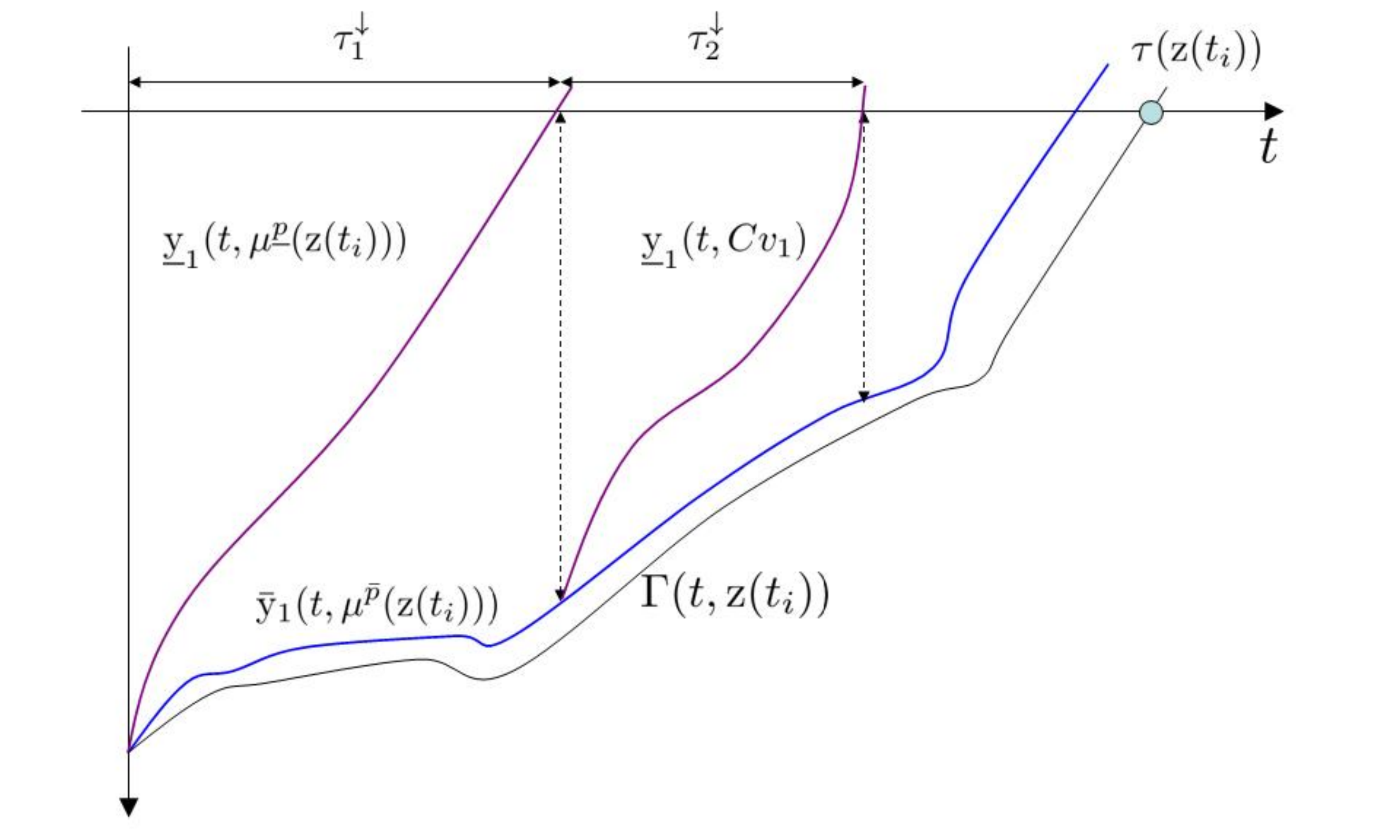}
   \end{center}
 \caption{Iterative algorithm for the computation of the inter-execution times.}
\label{fig_iter_alg}
\end{figure}
\setlength{\belowcaptionskip}{0pt} 

For clarity of exposition, we first describe the main ideas behind the algorithm. The procedure is sketched in Figure~\ref{fig_iter_alg}. Both $\underline{\y}_{1}$ and $\overline{\y}_{1}$ upper bound the evolution of the exact evolution of the triggering condition $\Gamma$, hence their roots lower bound the exact solution $\tau(\z(t_i))$. The first bound for the inter-execution time can be computed directly using the low order model:
\eq{
\label{first_t}
\tau^{\downarrow}_1 = \min\arg_t(\underline{\y}_1(t,\mu^{\underline{p}}(z(t_i))) = 0)
}
This value represents a bound on the inter-execution time $\tau(\z(t_i))$ defined by the triggering condition, $\Gamma(\z(\tau(\z(t_i)),\z(t_i)))=0$. Since $\underline{y}$ does not represent in general a tight bound for the evolution of $\Gamma$, it is expected that the triggering condition has not been satisfied at $\tau^{\downarrow}_1$. To refine the self-triggered time in~\eqref{first_t}, we evaluate the high order model at this instant of time:
\eq{
v_1 = \overline{\y}(\tau^{\downarrow}_1,\mu^{\overline{p}}(z(t_i)))
}
The high order model provides us with a lower bound for the gap between the exact inter-execution time as defined by the triggering condition and the first value obtained in~\eqref{first_t} through the low order model. We use this value $v_1$ as the new initial condition for $\underline{\y}_1$ to improve the initial guess $\tau^{\downarrow}_1$:
\eq{
\label{second_t}
\tau^{\downarrow}(\z(t_i)) = \tau^{\downarrow}_1 + \tau^{\downarrow}_2, \qquad \tau^{\downarrow}_2 = \min \arg_t(\underline{\y}_1(t,C v_1) = 0)
}
with $C = [\begin{array}{cc} I_{1\times \underline{p}} & 0_{1\times (\overline{p}-\underline{p})}\end{array}]$. The procedure can be repeated iteratively to obtain more accurate self-triggered times, as stated in the following theorem.
\begin{theorem}
\label{thm_iter}
Consider a dynamical system $\dot{\z}=Z(\z)$, homogeneous of degree $\xi>0$. The time $\tau^{\downarrow}(\z(t_i))$ given by: 
\eq{
\label{iter_sol_ideal}
\tau^{\downarrow}(\z(t_i)) = \sum_{l=1}^{n_{iter}} \tau_l^{\downarrow}
}
\eqs{
\label{iter_y}
\text{with} \qquad \tau_l^{\downarrow} &= &\min \arg_{\tau} \left(\underline{\y}_1(\tau, C v_{l-1}) =0\right), \qquad C = \left[\begin{array}{cc} I_{1\times \underline{p}} & 0_{1\times (\overline{p}-\underline{p})}\end{array}\right]\\
\label{def_v0}
\text{and} \qquad v_0 &=& \mu^{\overline{p}}(\z(t_i))\\
v_r &=& \overline{\y}(\tau^{\downarrow}_{r}, v_{r-1}), \quad r = 1, 2, \hdots, n_{\text{iter}}-1
}
represents for any $n_{\text{iter}}>0$ a lower bound for the inter-execution times \mbox{$\tau:\R^{2n} \rightarrow \R^+ \cup \{\infty\}$} implicitly given by the homogeneous triggering condition $\Gamma(\z(\tau(\z(t_i)),\z(t_i))) = 0$, with $\underline{y}$ and $\overline{y}$ as defined in~\eqref{def_ys}.
\end{theorem}

\begin{IEEEproof}
According to Lemma~\ref{lemma_approx_trig}, $\underline{\y}_1$ represents a bound for the evolution of the triggering condition:
\eq{
\label{eq_1}
\Gamma(\z(t,\z(t_i))) \leq \underline{\y}_1 (t,\mu^{\underline{p}}(\z(t_i))). 
}
Hence we can use $\underline{\y}_1$ to compute a first bound for the inter-execution time $\tau(\z(t_i))$:
\eq{
\label{eq_2}
\tau^{\downarrow}_1:=  \min\arg_t(\underline{\y}_1(t,\mu^{\underline{p}}(z(t_i))) = 0) \leq \tau(\z(t_i)):= \min\arg_t (\Gamma(\z(t,\z(t_i)))=0).
}
This bound will be refined using the high order model. According to Lemma~\ref{lemma_approx_trig}, $\overline{\y}$ also represents a bound for the evolution of the triggering condition and its Lie derivatives:
\eq{
\label{eq_3}
{\cal L}_Z^{l-1} \Gamma(\z(t, \z(t_i))) \leq \overline{\y}_l (t,\mu^{\overline{p}}(\z(t_i))),\quad l=1,\hdots,\overline{p}.
}
We define
\eq{
v_1 := \overline{\y} (\tau^{\downarrow}_1,\mu^{\overline{p}}(\z(t_i)))
}
which represents a bound for $\Gamma$ and its first $\overline{p}$ derivatives at the time instant $\tau^{\downarrow}_1$.
We shall use this bound as the new initial condition for $\underline{\y}_1$. We now rewrite inequality~\eqref{eq_2} using $v_1$:
\eq{
\label{eq_4}
\tau^{\downarrow}_2:=  \min\arg_t(\underline{\y}_1(t,Cv_1) = 0) \leq \tau_2 := \min\arg_t (\Gamma(\z(t,z_1))=0)
}
where $z_1$ represents the extended state after flowing along the vector field $Z$ for $\tau^{\downarrow}_1$ units of time:
\eq{
z_1 = \z(\tau^{\downarrow}_1,\z(t_i)).
}
Notice that inequality~\eqref{eq_4} holds since Lemma~\ref{lemma_approx_trig} also applies whenever the initial conditions satisfy:
\eq{
{\cal L}_Z^{l-1} \Gamma(z_0) \leq y_{0_l}, \quad l=1,\hdots,\overline{p}.
}
which is the case in~\eqref{eq_4} since:
\eq{
{\cal L}_Z^{l-1} \Gamma(z_1) \leq v_{1_l},\quad l=1,\hdots,\overline{p}.
}
Merging equations~\eqref{eq_2} and~\eqref{eq_4} we obtain a second bound for the inter-execution time:
\eq{
\tau^{\downarrow}_1 + \tau^{\downarrow}_2 \leq \tau^{\downarrow}_1 + \tau_2 = \tau(\z(t_i)).
}
Applying the procedure iteratively we finally obtain the final expression in~\eqref{iter_sol_ideal}, with $\tau_l^{\downarrow} = \min\arg_t(\underline{\y}_1(t,Cv_{l-1}) = 0)$.
\end{IEEEproof}
\begin{remark} 
There is no guarantee that $\overline{\y}<\underline{\y}$, for there does not exist a partial order on the approximations (that is, a high order approximation does not need to be more accurate than a low order approximation for all points in the extended state space). Hence, the bound does not need to increase monotonically with respect to the number of iterations. However, the iterative algorithm can be easily modified to account for such possibility: if $v_{i_1}>0$ at some step $i$, the computation should stop at this point, as the high order model estimates that the triggering condition has already been satisfied and thus will yield shorter times.
\end{remark}
In practice, since the expression for $\underline{y}$ is transcendental in $t$, a closed-form expression for $\tau^{\downarrow}_1$ and $\tau^{\downarrow}_2$ cannot be directly computed, as explained in Section~\ref{sec_self}. To find closed-form expressions for such times, we proceed as in Section~\ref{sec_self_trig_formula}: we use the scaling properties of the inter-execution times for a homogeneous system. The approximation of the isochronous manifold $\Omega$ using $\underline{y}_1$ is defined by:
\eq{
\underline{y}_1(t_*, \mu^{\underline{p}}(z_*))=0.
}
We find the intersection of the isochronous manifold with the homogeneous ray $z_* = \lambda \z(t_i)$:
\eq{
\label{eqq_1}
\underline{\y}_1(t_*,\mu^{\underline{p}}(\lambda z(t_i))) = \underline{\y}_1(t_*,\delta^{\lambda}(C o_0))= 0, \qquad o_0 = \mu^{\overline{p}}(\z(t_i))
}
since $\mu^{\underline{p}}(\lambda z(t_i)) = \delta^{\lambda}(\mu^{\underline{p}}( z(t_i)))$, according to Lemma~\ref{homog_lemma}. Equation~\eqref{eqq_1} is now polynomial of order $\underline{p}$ in $\lambda$ (after manipulating terms as in Section~\ref{sec_self_trig_formula}), and therefore a closed-form solution exists (since $\underline{p}\leq 3$). We denote such solution as $\lambda_1$:
\eq{
\label{eq_5}
\lambda_1 = \arg_{\lambda} (\underline{\y}_1(t_*,\delta^{\lambda}(C o_0)) = 0) \leq \arg_{\lambda} (\z(t_*, \lambda \z(t_i))=0).
}
Equation~\eqref{eq_5} implies that the point $\lambda_1 \z(t_i)$ has an inter-execution time bounded by $t_*$:
\eq{
\tau(\lambda_1 \z(t_i)) \geq t_*.
}
And using the scaling laws from Theorem~\ref{ScaleHom_gen} we compute the first bound for the inter-execution time at $\z(t_i)$:
\eq{
\tau(\z(t_i)) \geq \tau^{\downarrow} = \lambda_1^{\xi} t_*
}
which corresponds to equation~\eqref{final_formula}, with $q=\lambda_1^{\xi}$. Hence, for $n_{iter}=1$ we obtain the bound in Theorem~\ref{final_thm}. Likewise, we can use the same technique to compute the second iteration. As before, we define:
\eq{
\label{def_o}
o_1 = \overline{y}(t_*, \delta^{\lambda_1} (o_0))
}
which represents a bound for $\Gamma$ and its $\underline{p}$ Lie derivatives at the time instant $t_*$. We rewrite the left hand side of~\eqref{eqq_1} using $o_1$:
\eq{
\underline{\y}_1(t_*,\delta^{\lambda} (C o_1)) = 0.
}
We solve for $\lambda$ as before:
\eq{
\label{eq_10}
\lambda_2 = \arg_{\lambda} (\underline{\y}_1(t_*,\delta^{\lambda} (C o_1)) = 0)\leq \arg_{\lambda}(\z(t_*, \lambda z_a)=0)
}
where $z_a$ represents the extended state after flowing along the vector field $Z$ for $t_*$ units of time, starting with initial condition $\lambda_1 \z(t_i)$:
\eq{
z_a = \z(t_*,\lambda_1 \z(t_i)).
}
Equation~\eqref{eq_10} implies that the point $\lambda_2 z_a$ has an inter-execution time bounded by $t_*$:
\eq{
\tau(\lambda_2 z_a) \geq t_*.
}
Resorting to the scaling properties for $\tau$ again we can compute a bound for the inter-execution time at $\z(t_i)$:
\eq{
\tau(\lambda_2 z_a) = \lambda_2^{\xi} \tau(z_a) 
= \lambda_2^{\xi} \tau(\z(t_*,\lambda_1 \z(t_i))) 
=\lambda_2^{\xi} (\tau(\lambda_1 \z(t_i)) - t_*) = \lambda_2^{\xi} (\lambda_1^{\xi} \tau(\z(t_i)) - t_*) \geq t_*.
}
Finally, solving for $\tau(\z(t_i))$ we obtain:
\eq{
\tau(\z(t_i) \geq \lambda_1^{\xi} t_* + \lambda_2^{\xi} \lambda_1^{\xi} t_*.
}
Applying the procedure iteratively we obtain the following corollary, analogous to Theorem~\ref{thm_iter}:
\begin{corollary}
\label{coro_iter}
Consider a dynamical system $\dot{\z}=Z(\z)$, homogeneous of degree $\xi>0$. The time $\tau^{\downarrow}(\z(t_i))$ given by: 
\eq{
\label{iter_sol}
\tau^{\downarrow}(\z(t_i)) = \left(\sum_{l=1}^{n_{iter}} \left(\prod_{j=1}^l \lambda^{\xi}_j\right) \right) t_*
}
\eqs{
\label{iter_z}
\text{with} \qquad \lambda_j &=& \min \arg_{\lambda} \left(\underline{\y}_1(t_*, \delta^{\lambda} (C o_{j-1})) =0 \right), \\
\label{def_os}
\text{and} \qquad o_0 &=& \mu^{\overline{p}}(\z(t_i))\\
o_r &=& \overline{\y}(t_*,\delta^{\lambda_{r}} (o_{r-1})), \quad r = 1, 2, \hdots, n_{\text{iter}}-1.
}
represents for any $n_{\text{iter}}>0$ a lower bound for the inter-execution time \mbox{$\tau:\R^{2n} \rightarrow \R^+ \cup \{\infty\}$}  implicitly given by the homogeneous triggering condition $\Gamma(\z(\tau(\z(t_i)),\z(t_i))) = 0$, with $\underline{y}$ and $\overline{y}$ as defined in~\eqref{def_ys}.
\end{corollary}
\begin{remark}
The proposed algorithm has the property of being anytime, since it can be stopped at any iteration 
and still provide a bound for the inter-execution times. Other numerical algorithms, such as any Householder's method~\cite{householder1970numerical}, could be used to compute an approximation $\hat{\lambda}$ for the solution of the polynomial equation~\eqref{final_pol_eq}; then, the inequality $\overline{y}_1(\hat{\lambda})<0$ needs to be checked to guarantee that the solution indeed bounds the exact solution. However, in this approach no conclusions can be reached should the inequality not be satisfied. 
\end{remark}

\subsection{Example}

We consider now the control of a rigid body, previously studied in~\cite{tabuada10_tac}. The state space representation of such system with two inputs can be simplified to the form:
\begin{eqnarray}
\dot{\mathrm{x}}_1& = & \mathrm{u}_1 \nonumber \\
\dot{\mathrm{x}}_2& = & \mathrm{u}_2 \nonumber \\
\dot{\mathrm{x}}_3& = & \mathrm{x}_1 \mathrm{x}_2.
\end{eqnarray}
The following nonlinear feedback law was designed in~\cite{byrnes1989nra} to render the system globally asymptotically stable:
\begin{eqnarray}
u_1 &=& -x_1 x_2 - 2 x_2 x_3  -x_1 - x_3 \nonumber \\
u_2 &=& 2 x_1 x_2 x_3 + 3 x_3^2 -x_2.
\end{eqnarray}
We use the same triggering condition as in~\cite{tabuada10_tac}, given by:
\eq{
\label{trig_cond_rig_body}
|e|^2 = 0.79^2 \sigma^2 |x|^2,\qquad \sigma\in]0,1[.
}
Since the system is not homogeneous, we introduce an auxiliary variable $w$, as explained in Section~\ref{homog_subsec}. The homogenized closed loop system becomes:
\begin{equation*}
\label{example_homogenized}
\dot{\z}=Z(\z)=\left[\begin{array}{ccc} -(\mathrm{x}_1+\mathrm{e}_1) (\mathrm{x}_2+\mathrm{e}_2) \mathrm{w} - 2 (\mathrm{x}_2 + \mathrm{e}_2) (\mathrm{x}_3+\mathrm{e}_3)\mathrm{w}-(\mathrm{x}_1+\mathrm{e}_1)\mathrm{w}^2 - (\mathrm{x}_3+\mathrm{e}_3) \mathrm{w}^2 \\
2 (\mathrm{x}_1+\mathrm{e}_1) (\mathrm{x}_2+\mathrm{e}_2) (\mathrm{x}_3+\mathrm{e}_3) + 3 (\mathrm{x}_3+\mathrm{e}_3)^2\mathrm{w}- (\mathrm{x}_2 + \mathrm{e}_2)\mathrm{w}^2\\
\mathrm{x}_1 \mathrm{x}_2 \mathrm{w}\\
(\mathrm{x}_1+\mathrm{e}_1) (\mathrm{x}_2+\mathrm{e}_2) \mathrm{w} + 2 (\mathrm{x}_2 + \mathrm{e}_2) (\mathrm{x}_3+\mathrm{e}_3)\mathrm{w}+(\mathrm{x}_1+\mathrm{e}_1)\mathrm{w}^2 + (\mathrm{x}_3+\mathrm{e}_3) \mathrm{w}^2 \\
-2 (\mathrm{x}_1+\mathrm{e}_1) (\mathrm{x}_2+\mathrm{e}_2) (\mathrm{x}_3+\mathrm{e}_3) - 3 (\mathrm{x}_3+\mathrm{e}_3)^2\mathrm{w}+(\mathrm{x}_2 + \mathrm{e}_2)\mathrm{w}^2\\
-\mathrm{x}_1 \mathrm{x}_2 \mathrm{w}\\
0\end{array}\right]
\end{equation*}
with $\z=(\x,\e,\w)$ representing the extended state. It is not necessary to homogenize the triggering condition since it is already homogeneous.

We analyze the accuracy of the self-trigger times with respect to the number of iterations executed.
The degrees of the low order and high order models are $\underline{p}=3$ and $\overline{p}=4$ respectively.
Using this augmented vector field $Z$ and the triggering condition~\eqref{trig_cond_rig_body} we obtain a set of coefficients satisfying~\eqref{ineq_lies_lemma} for $p=\underline{p}$ and $p=\overline{p}$:
$$ p=\underline{p}: \quad \chi_0 = -73.2528,  \quad \chi_1= 1.7157,  \quad \chi_2 = 1.8299.$$
$$ p=\overline{p}: \quad  \chi_0 = -57.8151, \quad \chi_1= 1.4923,   \quad \chi_2 = 24.9920, \quad \chi_3 = 47.6313.$$
With this set of coefficients we can finally construct the polynomial equation in~\eqref{final_pol_eq_thm}, with $\z(t_i)=(\x(t_i),\e(t_i),1)$ (since $\w(t)\equiv1$).

The proposed self-triggering technique is evaluated for 25 initial conditions equally spaced along the boundary of a sphere of radius 1. Table~\ref{comparison_table_rigid_body} shows the inter-execution times (averaged over the considered initial conditions) when the number of iterations ranges from 1 to 3. It can be seen how the mismatch between the self-triggered times and the event-triggered times is quickly reduced as the number of iterations grows.

\begin{table}
\center
\begin{tabular}{c|ccc|c}
\multirow{2}{*} {$\sigma$} &  & self-triggered &  &  \multirow{2}{*}{event-triggered} \\ 
 &  $n_{\text{iter}}=1$ & $n_{\text{iter}}=2$ & $n_{\text{iter}}=3$ &\\
\hline\\
$0.50$ & 216.08 & 220.14 & 220.18 & 220.34\\
$0.65$ & 267.53 & 278.29 & 278.81 & 285.41\\
$0.80$ & 322.12 & 349.97 & 354.88 & 355.75 \\
\hline\\
\end{tabular}
\caption{Average time for the rigid body example (in ms.).}
\label{comparison_table_rigid_body}
\end{table}

\section{Selection of design parameters}







To conclude the paper, in this section we study the role played by the different design parameters of the proposed procedure. 

\subsection{Choice of the embedding degree}

In Section~\ref{sec_scal_laws} an embedding procedure was derived to render a system and a triggering condition homogeneous of any desired degree. In practice, one would like to choose the appropriate degree of the embedding to obtain accurate bounds for the inter-execution times. It was shown in Section~\ref{sec_homog_trig} that the scaling law does not depend on the degree of homogeneity of the triggering condition. 
In the following theorem we prove that embeddings of different degrees generate the same bounds for the inter-execution times of the original system.
\begin{theorem}
Let \mbox{$\tau_1:\R^{2n+1} \rightarrow \R^+ \cup \{\infty\}$}  be the inter-execution times of the system:
\eq{
\label{homogenized1st}
\left[\begin{array}{c} \dot{\z} \\ \dot{\w} \end{array}\right] = S_1(\z,\w)=\left[\begin{array}{c} \phantom{-}\w^{\xi_1+1} g(\w^{-1}\z) \\ 0 \end{array}\right]
}
implicitly defined by the homogeneous triggering condition $\Gamma(\w^{-1}\z	(\tau_1(z_0,w_0),z_0,w_0))=0$, and let \mbox{$\tau_2:\R^{2n+1} \rightarrow \R^+ \cup \{\infty\}$} be the inter-execution times for the system:
\eq{
\label{homogenized2nd}
\left[\begin{array}{c} \dot{\z} \\ \dot{\w} \end{array}\right] = S_2(\z,\w)=\left[\begin{array}{c} \phantom{-}\w^{\xi_2+1} g(\w^{-1}\z) \\ 0 \end{array}\right]
}
implicitly defined by $\Gamma(\w^{-1}\z(\tau_2(z_0,w_0),z_0,w_0))=0$. 
Consider two bounds $\tau^{\downarrow}_1$ and $\tau^{\downarrow}_2$ for $\tau_1$ and $\tau_2$, computed according to Theorem~\ref{final_thm} and valid in the sets $\Omega_1$ and $\Omega_2=\{(z,w)\in\R^{2n+1}\vert (w^{\xi_2/\xi_1-1}z,w^{\xi_2/\xi_1}) \in \Omega_1\}$ respectively. Then, the two bounds coincide for $w_0= 1$:
\eq{
\label{eq_bounds}
\tau_1^{\downarrow}(z_0,1) = \tau_2^{\downarrow}(z_0,1), \qquad \forall (z_0,1) \in \Omega_1 \cap \Omega_2.
}
\end{theorem}
\begin{IEEEproof}
Since the degree of homogeneity of the triggering condition $\Gamma$ does not play any role in the inter-execution times (as explained in Section~\ref{sec_homog_trig}), we assume it to be zero, without loss of generality. The bounds $\tau_1^{\downarrow}$ and $\tau_2^{\downarrow}$ are given by:
\eqs{
&\tau_j^{\downarrow} = q_j t_*, \qquad j=1,2&\\
&\sum_{i=0}^{p-1} \beta_{i_j}(z_0,1) q_j^i =0, \qquad \beta_{i_j}(z_0,1) = {\exp^{A_{p_j} t_*}}_{1(i+1)} \:({\mathcal L}_{S_j}^i \Gamma)(z_0,1), \quad i=0,\hdots,p-1. &\notag
}
Thus both bounds will coincide if $\beta_{i_1} = \beta_{i_2}$, that is, if $A_{p_1} =A_{p_2}$ and the Lie derivatives of both vector fields coincide for $w=1$, as imposed by the following equations:
\eq{
\label{statement1}
({\mathcal L}_{S_1}^i \Gamma)(z_0,1) = ({\mathcal L}_{S_2}^i \Gamma)(z_0,1) \qquad i=0,1,\ldots p-1\\
}
\eq{
\label{statement2}
\begin{split}
&({\mathcal L}_{S_1}^p \Gamma)(z_1,w_1) \leq \sum_{i=0}^{p-1} \chi_i ({\mathcal L}_{S_1}^i \Gamma)(z_1,w_1) \quad \forall (z_1,w_1) \in \Omega_1\\
\Leftrightarrow &({\mathcal L}_{S_2}^p \Gamma)(z_2,w_2) \leq \sum_{i=0}^{p-1} \chi_i ({\mathcal L}_{S_2}^i \Gamma)(z_2,w_2) \quad \forall (z_2,w_2) \in \Omega_2
\end{split}
}
where the coefficients $\chi_i$ are the same for both systems (in order to have the same $A_p$ matrix). 

To prove equality~\eqref{statement1}, we notice that the submanifold defined by $w=1$ is invariant under the flows of $S_1$ and $S_2$. Since $S_1$ and $S_2$ coincide for $w=1$, the evolution of the triggering condition in this submanifold is the same under both vector fields, and therefore the Lie derivatives coincide for any $(z,1)$.

To prove the double implication~\eqref{statement2}, we resort to the homogeneity properties of the Lie derivatives of $\Gamma$, as stated in Lemma~\ref{homog_lemma}. The proof relies on the fact that the vector fields coincide for $w=1$ and at the same time they are determined by their restriction to $w=1$.
To proceed, we represent a point $(z_1,w_1)\in\Omega_1$ as $\lambda_1(z_a, 1)$, and we rewrite the left hand side of~\eqref{statement2} as:
\eq{
\label{eqq_2}
{\mathcal L}_{S_1}^p \Gamma(\lambda_1 z_a, \lambda_1) \leq \sum_{i=0}^{p-1} \chi_i {\mathcal L}_{S_1}^i \Gamma(\lambda_1 z_a,\lambda_1)\qquad \forall \lambda_1>0 \vert (\lambda_1 z_a,\lambda_1)\in \Omega_1.
}
Since the $i$-th Lie derivative ${\mathcal L}_{S_1}^i \Gamma$ is homogeneous of degree $i \xi_1$, we can rewrite~\eqref{eqq_2} as:
\eq{
\lambda_1^{p \xi_1}  {\mathcal L}_{S_1}^p \Gamma(z_a, 1) \leq \sum_{i=0}^{p-1} \chi_i \lambda_1^{i \xi_1}  {\mathcal L}_{S_1}^i \Gamma(z_a,1) \qquad \forall \lambda_1>0 \vert (\lambda_1 z_a,\lambda_1)\in \Omega_1.
}
Defining $\lambda_2 := \lambda_1^{\xi_1/\xi_2}$ and using equality~\eqref{statement1} we obtain:
\eq{
\label{eqq_4}
\lambda_2^{p \xi_2}  {\mathcal L}_{S_2}^p \Gamma(z_a, 1) \leq \sum_{i=0}^{p-1} \chi_i \lambda_2^{i \xi_2}  {\mathcal L}_{S_2}^i \Gamma(z_a,1) \qquad \forall \lambda_2>0 \vert (\lambda_2^{\xi_2/\xi_1} z_a,\lambda_2^{\xi_2/\xi_1})\in \Omega_1.
}
Using the homogeneity properties again for the Lie derivatives equation~\eqref{eqq_4} becomes:
\eq{
\label{eqq_3}
{\mathcal L}_{S_2}^p \Gamma(\lambda_2 z_a, \lambda_2) \leq \sum_{i=0}^{p-1} \chi_i {\mathcal L}_{S_2}^i \Gamma(\lambda_2 z_a,\lambda_2)\qquad \forall \lambda_2>0 \vert (\lambda_2^{\xi_2/\xi_1} z_a,\lambda_2^{\xi_2/\xi_1})\in \Omega_1.
}
The point $(z_2,w_2):= (\lambda_2 z_a, \lambda_2)$ is related to $(z_1,w_1)$ according to:
\eq{
(z_2,w_2) = \lambda_2 \left(z_a,1\right) = \lambda_1^{\xi_1/\xi_2} \left(\frac{1}{\lambda_1}z_1,1\right) = \left(\lambda_1^{\xi_1/\xi_2-1} z_1,\lambda_1^{\xi_1/\xi_2}\right)=\left(w_1^{\xi_1/\xi_2-1} z_1,w_1^{\xi_1/\xi_2}\right).
}
Thus, inequality~\eqref{eqq_3} holds for all $(z_1,w_1)\in\Omega_1$, that is, for all $(w_2^{\xi_2/\xi_1-1}z_2,w_2^{\xi_2/\xi_1})\in\Omega_1$. Hence, we can conclude that the left hand side of~\eqref{statement2} for $(z_1,w_1)\in\Omega_1$ implies the right hand side of~\eqref{statement2} for $(z_2,w_2)\in\Omega_2$. The other side of the implication can be proved analogously by swapping the subindices $1$ and $2$. Therefore, $\tau^{\downarrow}_1$ and $\tau^{\downarrow}_2$ coincide for all those points $(z_0,1)$ where the bounds are valid, \ie, for all $(z_0,1)\in\Omega_1 \cap \Omega_2$
\end{IEEEproof}
The set $\Omega_2$ can be thought as an inflation of the set $\Omega_1$ as dictated by the different degrees of homogeneity. As both $\tau^{\downarrow}_1$ and $\tau^{\downarrow}_2$ coincide for all the points satisfying $w=1$, a bound for the inter-execution times of the dynamical system $\dot{\z}=Z(\z)$ under the triggering condition $\Gamma(\z(\tau,z_0))=0$ can be equally computed using either system~\eqref{homogenized1st} or~\eqref{homogenized2nd}. In other words, the degree of the homogenized system does not affect the accuracy of the bound for the inter-execution times.


\subsection{Choice of $t_*$}

In order to derive the self-triggering technique, we used the concept of isochronous manifolds in Section~\ref{iso_mnf_sect}. It was proven in Proposition~\ref{existence_iso_mnf} that such manifolds exist for \textit{any} value of $t_*$. This parameter has to be selected offline since computing $\exp^{A_p t_*}$ online would imply a high computational cost. Should the exact computation of this manifold be possible, any value of $t_*$ would yield the exact times. However, such exact computation is in general not feasible. Lemma~\ref{lemma_approx_trig} describes a way to compute one approximation bounding the exact isochronous manifold. Hence, a natural question is whether the choice of $t_*$ plays a role in the accuracy of the times, and, if so, which value of $t_*$ generates the largest times (and thus the most accurate bound).

To proceed, we analyze the role played by $t_*$ in equation~\eqref{final_pol_eq}:
\eq{
\label{eq1_dep_t}
\sum_{i=0}^{p-1} \beta_i(t_*,\z(t_i))  q^i = \sum_{i=0}^{p-1}{\exp^{A_p t_*}}_{1(i+1)} \:{\mathcal L}_Z^i \Gamma(\z(t_i)) q^i = 0.
}
Taking into consideration the special structure of the $A_p$ matrix (as defined in~\eqref{def_Ap}), we can write the Taylor expressions for the elements of $\exp^{A_p t_*}$ that appear in~\eqref{eq1_dep_t}:
\eq{
\label{eq2_dep_t}
{\exp^{A_p t_*}}_{1(i+1)} = \sum_{j=0}^{\infty} \frac{1}{j!} t_*^j (A_p^j)_{1(i+1)} = \frac{1}{i!} t_*^i + O(t_*^{p})
}
where $O(t_*^p)$ denotes terms that depend on $n$-powers of $t_*$, for $n\geq p$. Using~\eqref{eq2_dep_t}, equation~\eqref{eq1_dep_t} becomes: 
\eq{
\label{eq3_dep_t}
\sum_{i=0}^{p-1} \frac{1}{i!} t_*^i \:{\mathcal L}_Z^i \Gamma(\z(t_i)) q^i + \sum_{i=0}^{p-1}
O(t_*^{p}){\mathcal L}_Z^i \Gamma(\z(t_i)) q^i = 0.
}
Since we are interested in the dependence of $\tau^{\downarrow}$ on $t_*$, we rewrite equation~\eqref{eq3_dep_t} using $\tau^{\downarrow}=q t_*$: 
\eqs{
\label{eq4_dep_t}
&&\sum_{i=0}^{p-1} \frac{1}{i!} (\tau^{\downarrow})^i \:{\mathcal L}_Z^i \Gamma(\z(t_i)) + \sum_{i=0}^{p-1}
O(t_*^{p}){\mathcal L}_Z^i \Gamma(\z(t_i)) \left(\frac{\tau^{\downarrow}}{t_*}\right)^i = \notag\\
&&\sum_{i=0}^{p-1} \frac{1}{i!} (\tau^{\downarrow})^i \:{\mathcal L}_Z^i \Gamma(\z(t_i)) + \sum_{i=0}^{p-1}
O(t_*^{p-i}){\mathcal L}_Z^i \Gamma(\z(t_i)) (\tau^{\downarrow})^i = 0.
}
The first term of the previous summation corresponds to the first $p$ terms of the Taylor series of the triggering condition $\Gamma$. Therefore, the self-triggering formula that has been proposed in this paper represents a modified Taylor series that guarantees a lower bound for the inter-execution times. Likewise, we can obtain a similar expression for an upper bound $\tau^{\uparrow}$, as mentioned at the beginning of Section~\ref{sec_acc}:
\eq{
\label{eq5_dep_t}
\sum_{i=0}^{p-1} \frac{1}{i!} (\tau^{\uparrow})^i \:{\mathcal L}_Z^i \Gamma(\z(t_i)) + \sum_{i=0}^{p-1}
O(t_*^{p-i}){\mathcal L}_Z^i \Gamma(\z(t_i)) (\tau^{\uparrow})^i = 0.
}
Since the first terms of~\eqref{eq4_dep_t} and~\eqref{eq5_dep_t} coincide, the gap between the expressions for the upper bound and the lower bound is a function of $O(t_*)$, and therefore the mismatch between $\tau^{\uparrow}$ and $\tau^{\downarrow}$ increases with $t_*$. However, $t_*$ cannot be chosen as small as desired. The set $\Omega^{\downarrow}_{t_*}$ defined in~\eqref{def_iso_mnf} represents a bound for the isochronous manifold, valid in the set $M$ where the inequality~\eqref{ineq_lies_lemma} holds. As $t_*$ decreases, the diameter of the set $\Omega^{\downarrow}_{t_*}$ increases, since times enlarge as we move far from the origin (as dictated by the scaling law~\eqref{scal_eq_constant}). Hence, $t_*$ has to be chosen so that $\Omega^{\downarrow}_{t_*}$ is contained in the operating region $M$. Therefore, a sensible choice for $t_*$ is as follows:
\eq{
\label{t_star_choice}
t_* = \min \{t\in\R^+ \,\,\vert\,\, \Omega^{\downarrow}_{t}\subseteq  M\}.
}
Since the computation of~\eqref{t_star_choice} can be difficult in general, in practice $t_*$ can be chosen as the minimum value for which $\Omega^{\downarrow}_{t}\subseteq  M$ can be guaranteed. Notice that the choice of $t_*=0$ is not a possibility, since there does not exist a set of coefficients $\chi_i$ such that~\eqref{ineq_lies_lemma} holds globally. Indeed, the left hand side of the inequality~\eqref{ineq_lies_lemma} is a homogeneous function of higher degree than the homogeneous functions on the right hand side, thus there always exist points where ${\mathcal L}_{Z}^p \Gamma$ is greater than any linear combination of the previous Lie derivatives.

\section{Conclusions}
In this paper we developed a new self-triggering technique for nonlinear systems. Our previous work was extended in two directions: first, we enlarged the class of systems and triggering conditions that can be considered. 
Second, the new self-triggered formula reduces the conservativeness of previous techniques by resorting to the concept of isochronous manifolds, herein introduced. 
These results help consolidating self-triggered control as an interesting alternative to event-triggered control whenever the latter is not implementable. Moreover, even in those implementations where event-triggered control is preferred (\eg for robustness purposes), the self-triggering approach provides the tools to analyze the time requirements imposed by event-triggered control.

\bibliographystyle{alpha}
\bibliography{self_iso_mnf_bib}

\end{document}